\documentclass[review,hidelinks,onefignum,onetabnum,12pt]{article}
\usepackage{amsfonts}
\usepackage{graphicx}
\graphicspath{{figures/}}
\usepackage{color}      
\usepackage{hyperref}   
\usepackage{enumerate}
\usepackage{listings}
\usepackage{xcolor}
\usepackage{slashed}
\usepackage{amsmath}    
\usepackage{amssymb}
\usepackage{subfigure}
\usepackage{mathrsfs}
\usepackage{bbold}
\usepackage{tikz}
\usepackage[margin=1in,footskip=0.3in]{geometry}

\definecolor{mygreen}{rgb}{0,0.6,0}
\definecolor{mygray}{rgb}{0.5,0.5,0.5}
\definecolor{mymauve}{rgb}{0.58,0,0.82}

\newcommand{\im}{\mathbb{i}}

\title{Stability analysis of discontinuous Galerkin with \\ a high order embedded boundary treatment \\ for linear hyperbolic equations} 

\author{Mirco Ciallella\thanks{Laboratoire Jacques-Louis Lions, Universit\'e Paris Cit\'e, CNRS, UMR 7598, Paris, France 
  (\href{mailto:mirco.ciallella@u-paris.fr}{mirco.ciallella@u-paris.fr}).}
}
\date{}

\usepackage{amsopn}

\begin{document}

\maketitle

\begin{abstract}
Embedded, or immersed, approaches have the goal of reducing to the minimum 
the computational costs associated with the generation of body-fitted meshes by only employing fixed, 
possibly Cartesian, meshes over which complex boundaries can move freely. 
However, this boundary treatment introduces a geometrical error of the order of the mesh size that, 
if not treated properly, can spoil the global accuracy of a high order discretization, herein based on discontinuous Galerkin. 
The shifted boundary polynomial correction was proposed as a simplified version of the shifted boundary method, 
which is an embedded boundary treatment based on Taylor expansions to deal with unfitted boundaries. 
It is used to accordingly correct the boundary conditions imposed on a non-meshed boundary to compensate the aforementioned geometrical error, and reach high order accuracy. 
In this paper, the stability analysis of discontinuous Galerkin methods coupled with the shifted boundary polynomial correction 
is conducted in depth for the linear advection equation, by visualizing the eigenvalue spectrum of the high order discretized operators. 
The analysis considers a simplified one-dimensional setting by varying the degree of the polynomials and 
the distance between the real boundary and the closest mesh interface. 
The main result of the analysis shows that the considered high order embedded boundary treatment
introduces a limitation to the stability region of high order discontinuous Galerkin methods with explicit time integration, which becomes more and more important when using higher order methods.
The implicit time integration is also studied, showing that the implicit treatment of the boundary condition allows one to overcome such limitation and achieve an unconditionally stable high order embedded boundary treatment.
All theoretical results found by the analysis are validated through a large set of numerical experiments.

\end{abstract}



\section{Introduction}

This paper focuses on the stability analysis of high order discontinuous Galerkin methods \cite{reed1973triangular,cockburn2001runge} 
for the linear advection initial boundary value problem in one dimension:
\begin{equation}
    \label{eq:linear_advection_intro}
    \begin{cases}
    \partial_t u(x,t) + \partial_x u(x,t) = 0, \quad & x\in[x_L,x_R], \; t\geq 0 \\
    u(x_L,t) = u_D(t),\quad & t\geq 0, \\
    u(x,0) = u_0(x),\quad & x\in[x_L,x_R],
    \end{cases}
\end{equation}
where $u$ is the advected quantity and the advection velocity is set equal to one.
$x$ and $t$ represent the standard space and time variables, respectively.

The potential of high order methods in providing accurate solutions with reduced 
computational cost has been widely recognized \cite{wang2013high}.
However, ensuring consistency and stability in boundary conditions remains a challenge, and the difficulty
even increases when dealing with complex curved geometries.
As a matter of fact, accurately representing geometric features of the computational domain
can be as important as the numerical method itself \cite{bassi1997high}.
This is due to the fact that the global accuracy of the numerical solution is given by both
the accuracy of the numerical method and the accuracy of the geometry representation.

A widely used approach in high fidelity simulations is the isoparametric approach, which relies on 
high order polynomial reconstructions of the same degree as the numerical method to approximate the 
given geometry \cite{bassi1997high}. This provides optimal convergence rates and 
improved accuracy in the presence of curved boundaries, even when using coarse meshes.
Despite notable progress in this area, the generation and simulation of problems with 
high order curvilinear meshes remains a complex task: nonlinear mappings, complex quadrature rules,
and a higher number of required degrees of freedom are just a few of the challenges that need to be addressed;
additional difficulties are encountered in the mesh generation itself \cite{persson2009curved}, which needs
sophisticated techniques to ensure mesh validity and quality, and remains a multifaceted problem 
for realistic applications \cite{moxey2015isoparametric,coppeans2025aerodynamic}.

The goal of embedded boundary methods, sometimes also referred to as immersed or unfitted in different communities and/or numerical contexts, 
is to reduce to the minimum any effort related to mesh generation and mesh motion by only employing simple fixed (possibly Cartesian) meshes,
over which complex boundaries can move freely.
However, such boundary treatments introduce a geometrical error proportional to $h$, where $h$ is the mesh size.
By doing so, all challenges shift entirely from mesh treatment to the formulation of accurate and stable boundary conditions.
The philosophy behind these methods was first introduced in the pioneering work of Peskin 
\cite{peskin1977numerical,peskin2002immersed}, whose objective was to study the fluid-structure interaction inside the heart.
In the last fifty years, there has been a huge effort to develop and improve new embedded, immersed and unfitted methods 
with the same objective of cutting to the bone the computational costs imposed by mesh generation and motion when dealing with complex geometries.
Some of these works can be found in \cite{mittal2005immersed,fadlun2000combined,hansbo2002unfitted,fedkiw1999non,berger2012progress,monasse2012conservative,coco2013finite}, 
but unfortunately it is impossible to mention them all.
Since most of these methods are only first or second order accurate, 
the growing interest in high order methods has brought researchers to further investigate the possibility of retaining
high order accuracy on boundaries, even when the geometry is not meshed.  
Some recent works in this direction can be found in \cite{tan2010inverse,atallah2022high,giuliani2022two}. 

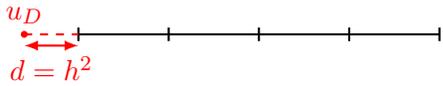
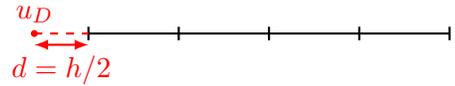
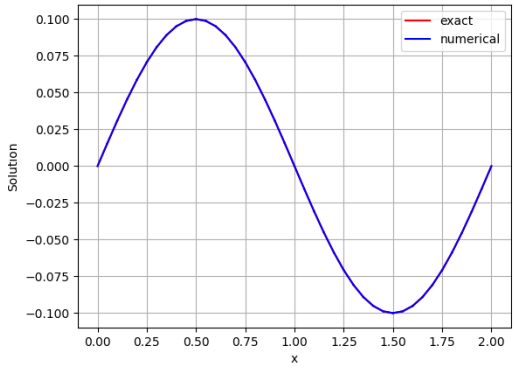
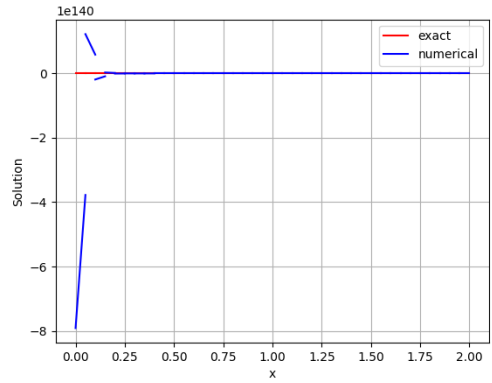
\begin{figure}
\centering
\subfigure[Geometrical error $h^2$]{
\begin{tikzpicture}[>=latex,scale=0.6]
    \draw[thick] (0,-0.15) -- (0,0.15);
    \draw[thick] (2,-0.15) -- (2,0.15);
    \draw[thick] (4,-0.15) -- (4,0.15);
    \draw[thick] (6,-0.15) -- (6,0.15);
    \draw[thick] (8,-0.15) -- (8,0.15);
    \draw[thick] (0,0) -- (2,0);
    \draw[thick] (2,0) -- (4,0);
    \draw[thick] (4,0) -- (6,0);
    \draw[thick] (8,0) -- (6,0);
    \draw[thick,red,dashed] (-1.2,0) -- (0,0);
    \draw[thick,fill=red,red] (-1.2,0) circle (1.5pt);
    \draw[thick,<->,red] (-1.2,-0.25) -- (0,-0.25);
    \node[above,red,font=\small] at (-1.2,0) {$u_D$};
    \node[below,red,font=\small] at (-0.6,-0.25) {$d=h^2$};
\end{tikzpicture}
}\hfill
\subfigure[Geometrical error $h/2$]{
\begin{tikzpicture}[>=latex,scale=0.6]
    \draw[thick] (0,-0.15) -- (0,0.15);
    \draw[thick] (2,-0.15) -- (2,0.15);
    \draw[thick] (4,-0.15) -- (4,0.15);
    \draw[thick] (6,-0.15) -- (6,0.15);
    \draw[thick] (8,-0.15) -- (8,0.15);
    \draw[thick] (0,0) -- (2,0);
    \draw[thick] (2,0) -- (4,0);
    \draw[thick] (4,0) -- (6,0);
    \draw[thick] (8,0) -- (6,0);
    \draw[thick,red,dashed] (-1.2,0) -- (0,0);
    \draw[thick,fill=red,red] (-1.2,0) circle (1.5pt);
    \draw[thick,<->,red] (-1.2,-0.25) -- (0,-0.25);
    \node[above,red,font=\small] at (-1.2,0) {$u_D$};
    \node[below,red,font=\small] at (-0.6,-0.25) {$d=h/2$};
\end{tikzpicture}
}
\subfigure[Simulation for $d=h^2$]{\includegraphics[width=0.45\textwidth]{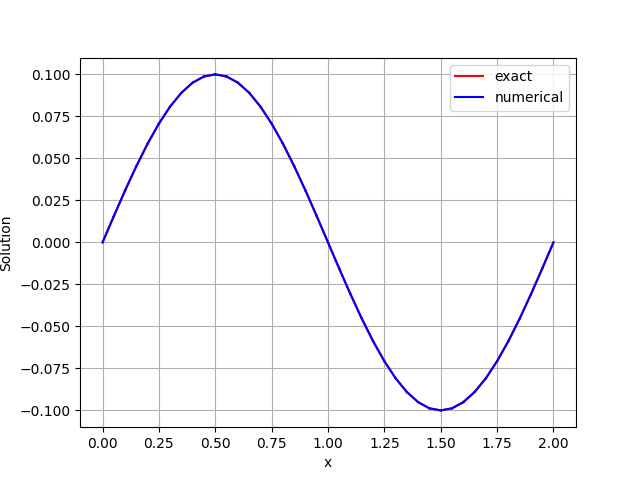} \label{fig:1c}} \hfill
\subfigure[Simulation for $d=h/2$]{\includegraphics[width=0.45\textwidth]{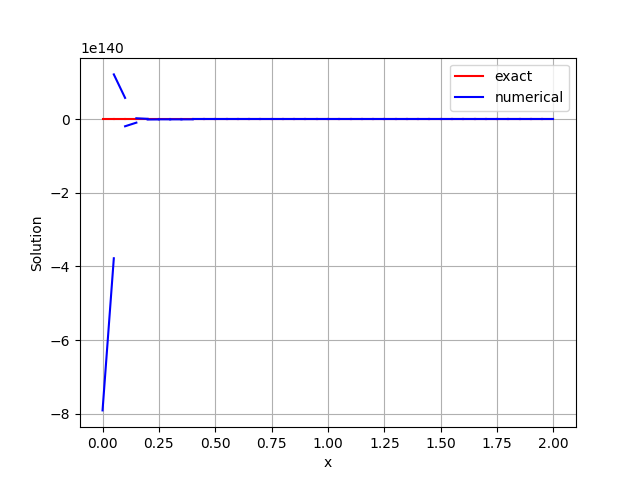}\label{fig:1d}}
\caption{These numerical results, obtained using a third order accurate $\mathbb{P}^2$ DG approximation of problem \eqref{eq:linear_advection_intro}, show the impact of the distance on the extrapolation of Dirichlet boundary condition using the shifted boundary embedded approach. On the left, the distance $d$ between the real boundary and the mesh is taken equal to $h^2$, mimicking a curved boundary discretized with simplicial elements, and the solution stays stable. On the right, $d$ is taken equal to $h/2$, mimicking a fully non-meshed configuration, and the solution is unstable.}
\label{fig:introexperiment}
\end{figure}

This paper contributes to the development and understanding of high-order embedded boundary methods applied to hyperbolic partial differential equations (PDEs). 
In particular, this work proposes to study in detail the stability analysis of discontinuous Galerkin (DG) methods coupled with 
the shifted boundary polynomial correction \cite{ciallella2023shifted} for problems with fully non-meshed boundary configurations, 
meaning when the distance between the exact boundary and the closest mesh edge is of the order $h$.
The shifted boundary polynomial correction was introduced in \cite{ciallella2023shifted},
as an efficient alternative to the high order shifted boundary method \cite{atallah2022high}.
Some applications of the shifted boundary method can be found in 
\cite{main2018shifted,atallah2022high,ciallella2020extrapolated,ciallella2022extrapolated,assonitis2022extrapolated}.
The shifted boundary polynomial correction was developed to achieve high order boundary conditions for hyperbolic problems with curved boundaries discretized by simplicial meshes (triangles or tetrahedrons),
meaning  when the distance between the exact boundary and the closest mesh element is of the order $h^2$.
Other methods to achieve high order consistency on simplicial elements have been presented in
\cite{costa2018very,burman2018cut,ciallella2023shifted,ciallella2024very}.
In \cite{ciallella2023shifted}, it was shown that it is possible to rewrite the corrected shifted boundary condition as a simple polynomial correction, while avoiding the explicit computation of high order derivative terms. 
More recently, it has also been developed for moving curved boundaries discretized with arbitrary-Lagrangian-Eulerian finite volume schemes on simplicial meshes \cite{BOSCHERI2025114215},
showing to be robust once again in several configurations involving a geometrical error of order $h^2$.
This approach, as the standard shifted boundary formulation, has also been shown to be valid in fully non-meshed configurations for high order finite element approximations of elliptic PDEs \cite{visbech2025spectral,atallah2022high}.
Contrary to the elliptic case, hyperbolic PDEs present additional challenges related to the complete absence of viscosity in the model, which makes numerical methods more prone to develop instabilities.
With a simple numerical experiment in one dimension, it is possible to show how the shifted boundary extrapolation affects the stability of the numerical method when varying the distance between the real boundary and the closest mesh interface.
In figure \ref{fig:introexperiment}, the numerical solution of  a third order accurate $\mathbb{P}^2$ DG approximation of problem \eqref{eq:linear_advection_intro} is presented with two configurations. Figure \ref{fig:1c} considers that the geometrical error is taken equal to $h^2$, mimicking a curved boundary discretized with simplicial elements, and shows that the solution is stable and converges to the exact one. Instead, Figure \ref{fig:1d} shows that instabilities arise when taking a fully non-meshed configuration, with a geometrical error equal to $h/2$.
This experiment motivates the need to carry out a thorough stability analysis.
Herein, the goal is to study the eigenvalue spectrum of the discretized operators of high order DG approximations of linear hyperbolic equations, 
when introducing the shifted boundary polynomial correction to design the corresponding high order embedded boundary condition.
The approach of eigenspectrum visualization has shown to be valid in different contexts \cite{M2AN_2015__49_1_39_0,michel2021spectral,yang2025inverse}, especially for high order methods whose analytical study may be extremely difficult to perform.  
In particular, it is the first time in the literature that such an analysis is proposed for the shifted boundary embedded treatment with high order DG methods in the context of hyperbolic problems. 
Works on similar subjects, found in the literature of hyperbolic PDEs, include the analysis of the inverse Lax-Wendroff boundary treatment for high order finite difference schemes \cite{li2016stability,M2AN_2015__49_1_39_0}. 
The stability analysis of the inverse Lax-Wendroff boundary treatment for DG methods was proposed very recently in \cite{yang2025inverse}.
For low order methods, a stability analysis with a cut-cell boundary treatment was presented in \cite{BERGER2015180}.

The main result of the analysis presented herein shows that the considered high order embedded boundary treatment
introduces a notable limitation to the stability region of high order DG methods with explicit time integration, 
and that this limitation becomes more and more important when using higher order methods.
The analysis also shows that the implicit time integration allows one to overcome such constraints
and achieve an unconditionally stable high order embedded boundary treatment.
All results obtained from the proposed stability analysis are validated through a large set of numerical simulations where
several boundary configurations are considered, with accuracy varying between second and fourth order.

The rest of the paper is organized as follows.
In section \ref{sec:DGadvection} the high order DG discretization for the linear advection equation is presented.
Section \ref{sec:SBcorrection} describes the high order embedded boundary treatment considered herein, namely the shifted boundary polynomial correction.
The stability analysis are presented in section \ref{sec:stability}.
The analysis is applied first to DG with standard periodic boundary conditions in section \ref{sec:periodic}. 
Then, the DG methods with the shifted boundary polynomial correction and explicit time integration is studied in section \ref{sec:explicitSB},
while the implicit treatment of the high order embedded boundary is considered in section \ref{sec:implicitSB}.
Numerical tests that validate all the theoretical results are presented in section \ref{sec:tests},
and some conclusions are drawn in section \ref{sec:conclusions}.

\section{High order discontinuous Galerkin for linear advection equations}\label{sec:DGadvection}

The one-dimensional computational domain $\Omega$ is discretized into $N_e$ non-overlapping elements $\Omega_e$, such that $\Omega=\bigcup_{e=1}^{N_e} \Omega_e$.
Considering a one-dimensional uniform mesh, the cell $\Omega_e=[x_{e-1/2},x_{e+1/2}]$ has a fixed size of $\Delta x$. 
The numerical solution $u$ is approximated by $u_h$, which belongs to the global discrete space $V_h$ defined as:
\begin{equation}
V_h:= \bigoplus_{e=1}^{N_e} V_e, \quad \text{where} \quad V_e := \mathbb{P}^p(\Omega_e).
\end{equation}
Here, $V_e$ is the local approximation space on element $\Omega_e$, given by the space of polynomials of degree $p \ge 0$. Consequently, the numerical solution is a piece-wise polynomial that is discontinuous across element interfaces.

Within each element $\Omega_e$, the approximation $u_h$ is expressed as a linear combination of basis functions $\{\phi_j\}_{j=0}^{p}$ of the space $V_e$:
\begin{equation}
	\label{eq:discontinuous_approximation}
    u_h(x, t)|_e = \sum_{i=0}^{p} \phi^e_j(x) u_{e,j}(t) ,
\end{equation}
where $u_{e,j}(t)$ are the time-dependent degrees of freedom.
Herein, a modal DG formulation is considered, but the same analysis can be extended to nodal DG methods without loss of generality.

To study the impact of Dirichlet boundary conditions, the linear advection equation with a source term will be considered for the space discretization:
\begin{equation}
    \label{eq:linear_advection_source}
    \partial_t u(x,t) + \partial_x u(x,t) = s(x).
\end{equation}
In the following, the notation $(x,t)$ will be dropped to be concise. 

The elemental semi-discrete DG weak formulation is written by projecting \eqref{eq:linear_advection_source} 
and integrating by parts, 
\begin{equation}
	\label{eq:weak_formulation_integrated}
	\int_{\Omega_e} \phi^e_i \frac{\text{d} u_h}{\text{d} t}  \,\text{d}x - \int_{\Omega_e} \partial_x  \phi^e_i u_h \,\text{d}x + 
    \int_{\partial \Omega_e} \phi^e_i \hat u(u_h^-,u_h^+) n \,\text{d}\Sigma = \int_{\Omega_e} \phi^e_i s(x) \,\text{d}x , \quad\forall \phi^e_i \in V_e,
\end{equation}
where $n$ is the outward normal, and $\hat u(u_h^-,u_h^+)$ is a numerical flux that depends on the left and right reconstructed interface values.
In one-dimension, equation \eqref{eq:weak_formulation_integrated} can be simply rewritten as
\begin{equation}
	\label{eq:weak_formulation_flux}
	\int_{\Omega_e} \phi^e_i \frac{\text{d} u_h}{\text{d} t}  \,\text{d}x - \int_{\Omega_e} \partial_x \phi^e_i u_h \,\text{d}x + 
     \bigl[\phi^e_i \hat u(u_h^-,u_h^+)\bigr]^{e+1/2}_{e-1/2} = \int_{\Omega_e} \phi^e_i s(x) \,\text{d}x , \quad\forall \phi^e_i \in V_e.
\end{equation}
The numerical flux at the element interfaces is computed using a simple upwind Riemann solver $\hat u(u_h^-,u_h^+) = u_h^-$, 
since the advection speed is positive.

When replacing \eqref{eq:discontinuous_approximation} into equation \eqref{eq:weak_formulation_flux} for the linear problem \eqref{eq:linear_advection_source}, 
the semi-discrete problem for a generic element $\Omega_e$ can be written in matrix form as follows 
\begin{equation}
	\label{eq:discontinuous_galerkin_matrix_form}
	M_e  \frac{\text{d} \mathbf{u}_e}{\text{d} t} - K^s_e \mathbf{u}_e + K^R_e \mathbf{u}_e - K^L_e \mathbf{u}_{e-1} = \mathbf{s}_e, 
\end{equation}
where  $\mathbf{u}_e$ are the modes belonging to element $\Omega_e$.
The mass and stiffness matrices are defined as,
\begin{equation}
    \label{eq:mass_stiff_mat}
    M_e = \int_{\Omega_e} \phi^e_i \phi^e_j \,\text{d}x ,\qquad K^s_e = \int_{\Omega_e} \partial_x \phi^e_i \phi^e_j \,\text{d}x ,\quad\forall \phi^e_i, \phi^e_j \in V_e, 
\end{equation}
while interface matrices are defined as,
\begin{equation}
    \label{eq:interface_mat}
    K^R_e =  \phi_i^e(x_{e+1/2}) \phi_j^e(x_{e+1/2})  ,\qquad K^L_e = \phi_i^e(x_{e-1/2}) \phi_j^{e-1}(x_{e+1/2})  
\end{equation}
where $\phi_i^e(x_{e+1/2})$ represents the $i$-th basis of the element $\Omega_e$ evaluated in $x_{e+1/2}$.

\section{Shifted boundary polynomial correction} \label{sec:SBcorrection}

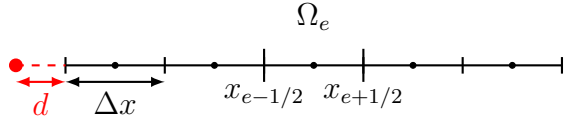
\begin{figure}
\centering
    \begin{tikzpicture}[>=latex,scale=1.1]

    \def\n{5} 
    \def\dx{1.2} 

    \foreach \i in {0,...,\n}{
    \draw[thick] (\i*\dx,-0.1) -- (\i*\dx,0.1);
    }

    \draw[thick] (0,0) -- (\n*\dx,0);
    \foreach \i in {1,...,\n}{
    \filldraw[black] ({(\i-0.5)*\dx},0) circle (1pt);
    }
    \draw[thick] (3*\dx,-0.2) -- (3*\dx,0.2);
    \draw[thick] (2*\dx,-0.2) -- (2*\dx,0.2);

    \node[below] at (3*\dx,-0.1) {$x_{e+1/2}$};
    \node[below] at (2*\dx,-0.1) {$x_{e-1/2}$};

    \node[above] at ({(3-0.5)*\dx},0.3) {$\Omega_e$};

    \draw[thick,red,dashed] (-0.6,0) -- (0,0);
    \draw[thick,fill=red,red] (-0.6,0) circle (2pt);
    \draw[thick,<->,red] (-0.6,-0.2) -- (0,-0.2);
    \node[below,red] at (-0.3,-0.2) {$d$};
    \draw[thick,<->] (0,-0.2) -- (1.2,-0.2);
    \node[below] at (0.6,-0.2) {$\Delta x$};

    \end{tikzpicture}
    \caption{Mesh configuration: tessellation of one-dimensional elements with a boundary condition (red circle) imposed at a distance $d$ from the left-most interface of the domain, which is considered as the surrogate boundary of the problem.}\label{fig:1dmesh}
\end{figure}

The aim of this work is to explore the effect of the shifted boundary correction, a high order embedded boundary treatment, when the distance $d$
from the real boundary to the left-most interface of the computational mesh is of the order $\Delta x$. This is a classical configuration that can occur 
when dealing with numerical simulations with boundaries treated as {\it immersed} starting from a Cartesian mesh, as can be also noticed in figure \ref{fig:1dmesh}.

\begin{figure}
\centering
    \begin{tikzpicture}[>=latex,scale=1.1]

    \draw[thick] (0,-0.2) -- (0,0.2);
    \draw[thick] (3,-0.2) -- (3,0.2);
    \draw[thick] (0,0.) -- (3,0);
    \draw[thick,red,dashed] (-1.6,0) -- (0,0);
    \draw[thick,fill=red,red] (-1.6,0) circle (2pt);
    \draw[thick,<->,red] (-1.6,-0.2) -- (0,-0.2);
    \node[below,red] at (-0.8,-0.2) {$-|d|$};
    \node[above,red,scale=1.2] at (-1.6,0.2) {$\bar x$};
    \node[above,scale=1.2] at (0,0.2) {$\tilde x$};

    \def\y{2} 

    \draw[thick] (0,-0.2-\y) -- (0,0.2-\y);
    \draw[thick] (3,-0.2-\y) -- (3,0.2-\y);
    \draw[thick] (0,0.-\y) -- (3,0-\y);
    \draw[thick,red,dashed] (+1.6,0-\y) -- (0,0-\y);
    \draw[thick,fill=red,red] (+1.6,0-\y) circle (2pt);
    \draw[thick,<->,red] (+1.6,-0.2-\y) -- (0,-0.2-\y);
    \node[below,red] at (+0.8,-0.2-\y) {$|d|$};
    \node[above,red,scale=1.2] at (1.6,0.2-\y) {$\bar x$};
    \node[above,scale=1.2] at (0,0.2-\y) {$\tilde x$};

    \end{tikzpicture}
    \caption{Boundary position: the real boundary position can be either outside of the first mesh element (above), or inside the first mesh element (below).
    In the first case, the distance will be taken as negative $d<0$, while for the second one the distance is taken positive $d>0$.}\label{fig:boundary}
\end{figure}
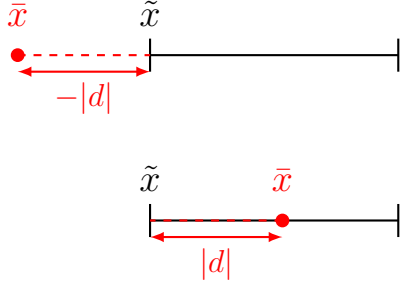

For multi-dimensional problems the mapping between the real boundary and the so-called {\it surrogate} one, which represents the union of the 
closest mesh edges to the real boundary, is not unique and can be prescribed using different approaches: closest-point, level-set, projections.
However, in one-dimensional configurations, the mapping between the two is uniquely identified, as shown in figure \ref{fig:boundary}. 
This analysis leaves the freedom to the real boundary, to be treated as an external or internal point of the mesh: meaning that the boundary
point could be either outside the first mesh element or inside it (see figure \ref{fig:boundary} for details). 
Consider $\bar x$ the position of the real boundary, while $\tilde x$ the position of the surrogate boundary. 
In this simplified configuration, the mapping between the two is simply $\bar x = \tilde x + \tilde d$, where $\tilde d$ is the distance vector.

Considering $n_e$ the unit outward normal of the left interface of the boundary element, the distance in one-dimension can be defined as a scalar quantity:
\begin{equation}\label{eq:distancedef}
d = - \langle n_e, \bar{x} - \tilde{x} \rangle,
\end{equation}
where $\langle \cdot,\cdot\rangle$ is the scalar product.
From now on, the vector notation $\tilde d$ will be dropped for simplicity, and the distance will be simply referred to as $d$.
Following this notation, the distance is allowed to vary within the range $d\in[-\Delta x, \Delta x]$.

In this section, the main notation of the shifted boundary polynomial correction to impose extrapolated boundary conditions is recalled.
The idea of the considered boundary treatment consists in imposing boundary conditions on the surrogate boundary rather than the real one, 
while retaining the high order accuracy thanks to a truncated Taylor expansion.
Hence, let $u_D$ be the prescribed Dirichlet boundary condition of the advected variable on the real boundary $\bar x$.
A modified boundary condition $u^\star$ for the surrogate boundary $\tilde x$ can be developed starting from the following Taylor expansion in one dimension:
\begin{equation}\label{eq:taylor}
    u(\bar x)=u(\tilde x + d) = u(\tilde x) + \sum_{m=1}^p \partial_x^{(m)} u(\tilde x) \frac{d^m}{m!}.
\end{equation}
Thus, a modified surrogate boundary condition can be defined with all the additional corrective terms and 
by imposing the real boundary condition $u(\bar x) = u_D$:
\begin{equation}
    u^\star(\tilde x) = u_D - \sum_{m=1}^p \partial_x^{(m)} u(\tilde x) \frac{d^m}{m!}.  
\end{equation}
The shifted boundary polynomial correction can be extracted naturally from equation \eqref{eq:taylor},
meaning that the whole set of high order correction terms can be recovered completely by simply evaluating the 
polynomial of the boundary cell in the point $\bar x = \tilde x+d$ and by making the difference with the solution evaluated in $\tilde x$.
Finally the modified Dirichlet boundary condition $u^\star$ on the surrogate boundary can be reformulated under the following correction:
\begin{equation}\label{eq:correction}
    u^\star(\tilde x) = u_D - (u(\tilde x + d) - u(\tilde x)) = u(\tilde x) - (u(\bar x) - u_D),
\end{equation}
which consists in shifting the surrogate boundary value with the polynomial error on the true boundary.
The main advantage of this formulation is its simplicity in the implementation, and the fact that 
all those corrective terms based on the solution derivatives are replaced by a simple evaluation of the solution on the real boundary.
The boundary condition at $x_L$ is then imposed at left-interface of the boundary cell through the numerical flux $\hat u (u^\star,u^+_h)= u^\star$, for an upwind flux. 
In particular, the evolution equation \eqref{eq:weak_formulation_flux} with an upwind numerical flux for the boundary element $[x_{1/2},x_{3/2}]$, which is the first mesh element, reads
\begin{equation}
\begin{split}
	\label{eq:weak_formulation_boundaryflux}
	&\int_{\Omega_1} \phi^1_i \frac{\text{d} u_h}{\text{d} t}  \,\text{d}x - \int_{\Omega_1} \partial_x \phi^1_i u_h \,\text{d}x + \\
    &\qquad\qquad \bigl[\phi^1_i(x_{3/2}) u_h(x_{3/2})|_1 - \phi^1_i(x_{1/2}) u^\star(x_{1/2}) \bigr] = \int_{\Omega_1} \phi^1_i s(x) \,\text{d}x . 
\end{split}
\end{equation}

In the following section, the analysis of the method will be performed by looking at the eigenvalue spectrum of the discretized operators,
and the stability region will be computed numerically by varying the CFL condition, $\text{CFL}=\Delta t/\Delta x$ in this simplified case, 
and the distance $d$.

\section{Stability analysis}\label{sec:stability}

To simplify the writing of this section, the analysis will be performed with a system made of two elements containing $p+1$ degrees of freedom each.
In contrast to the Von Neumann analysis, which assumes an infinite number of cells, the presence of boundary conditions makes the analysis dependent on the finite number of cells in the domain. Nevertheless, the qualitative behavior of the results remains largely unaffected by the specific number of cells.
This is also typical of other analysis involving boundary conditions that have been presented in the literature \cite{M2AN_2015__49_1_39_0}.
For explicit time integration schemes, the reference maximum CFL number will be that of the internal DG scheme, computed assuming a periodic domain. For a small number of cells the boundary conditions may locally enlarge the stability region compared to this reference value. However, when a sufficiently large number of cells is considered (as in practical cases) the overall stability region remains bounded by that of the internal scheme.

In section \ref{sec:periodic}, the first analysis is performed by enforcing standard periodic boundary conditions 
to recover the classical CFL stability constraints of DG schemes. 
This will prove the validity of the approach that will be used afterwards for the stability of high order DG methods
with the shifted boundary polynomial correction in section \ref{sec:explicitSB} (with explicit time integration) 
and in section \ref{sec:implicitSB} (with implicit time integration). 

To trace the stability region of the methods from the eigenvalue spectrum, 
it is necessary to write the system in a global matrix form of the kind
\begin{equation}\label{eq:semidiscrete} 
    \mathbf{M}\frac{\text{d} \mathbf{U}}{\text{d} t}  = \mathbf{K} \mathbf{U}, 
\end{equation} 
where $\mathbf{M}$ and $\mathbf{K}$ represent the global mass and stiffness matrices, respectively.
$\mathbf{U}$ is the global vector of modal degrees of freedom.

 To derive the eigenvalue problem, the explicit time integration of the semi-discrete system 
\eqref{eq:semidiscrete}, for a Runge-Kutta(4,4) gives 
\begin{equation}\label{eq:amplificationRK} 
    \mathbf{U}^{n+1} = \left( 1 + \mu + \frac{\mu^2}{2} + \frac{\mu^3}{6} +  \frac{\mu^4}{24}  \right)\mathbf{U}^n, 
\end{equation}
where $\mathbf{U}^n = \mathbf{U}(t^n)$ and $\mu=\lambda \Delta t$, with $\Delta t$ the time step and 
$\lambda$ the eigenvalue of semi-discrete operator $\mathbf{M}^{-1}\mathbf{K}$.
In this work the polynomial degree for the space discretization is taken up to $p=3$ and it is coupled with
a classical Runge-Kutta($p+1$,$p+1$) time integration scheme. 
Therefore, the stability region of the considered fully discrete problem is
$|z(\mu)|\leq 1$, where $z(\mu) = \sum_{k=0}^{p+1} \frac{\mu^k}{k!}$.

As for implicit integration, only implicit Euler will be considered to show in practice the beneficial effect of implicit time integration.

\subsection{Stability analysis: standard periodic boundary conditions}\label{sec:periodic}

\begin{figure}
    \centering
\begin{tikzpicture}[>=latex,scale=1.1]

\draw[thick,dashed] (0,0) -- (2,0);
\draw[thick] (2,0) -- (4,0);
\draw[thick] (4,0) -- (6,0);

\draw[thick] (0,-0.2) -- (0,0.2);

\foreach \i/\x/\label in { 1/2/{$x_{1/2}$}, 2/4/{$x_{3/2}$}, 3/6/{$x_{5/2}$}}
{
    \draw[thick] (\x,-0.2) -- (\x,0.2);
    \node[above] at (\x,0.25) {\label};
}

\node[below] at (1, -0.4) {$K_2$};
\node[below] at (3, -0.4) {$K_1$};
\node[below] at (5, -0.4) {$K_2$};

\end{tikzpicture}
\caption{Periodic boundary configuration used in the analysis of section \ref{sec:periodic}.}\label{fig:periodicDisegno}
\end{figure}
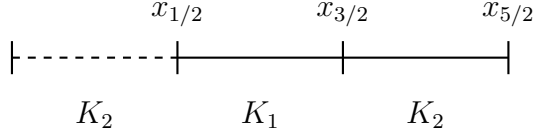

\begin{figure}[h!]
\centering
\subfigure{\includegraphics[width=0.7\textwidth]{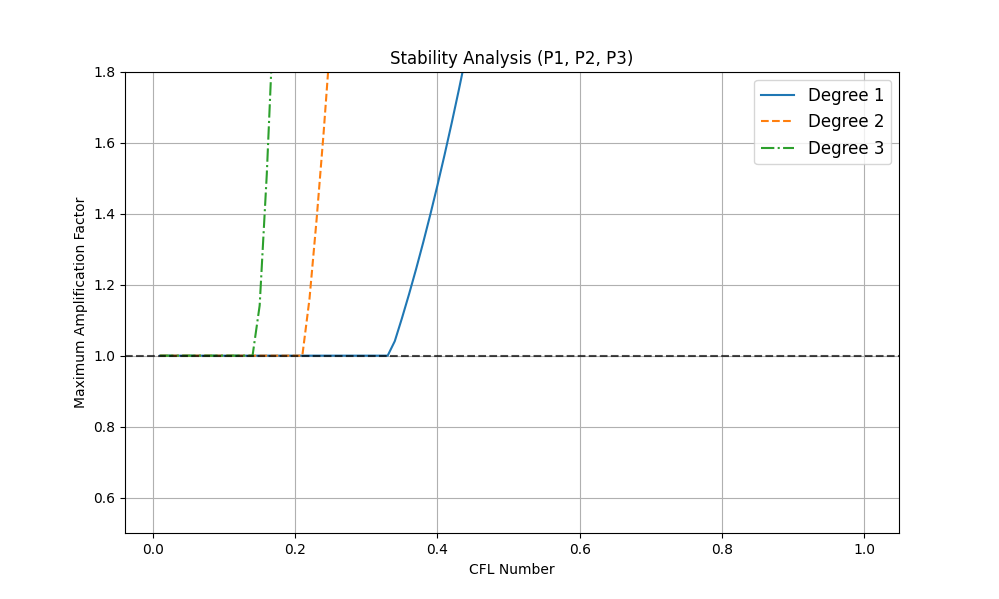}}
\caption{Amplification factor: periodic boundary  conditions  for the explicit discontinuous Galerkin method.}
\label{fig:Amplification_periodic_explicit}
\end{figure}

A system of two cells with periodic boundary conditions imposed at the left (of the first cell) and at the right (of the second cell) reads
\begin{equation}\label{eq:system_periodic}
    \begin{split}
	M_1 \frac{\text{d} \mathbf{u}_1}{\text{d} t} &= K^s_1 \mathbf{u}_1 - ( K^R_1 \mathbf{u}_1 - K^L_1 \mathbf{u}_2 ), \\
	M_2 \frac{\text{d} \mathbf{u}_2}{\text{d} t} &= K^s_2 \mathbf{u}_2 - ( K^R_2 \mathbf{u}_2 - K^L_2 \mathbf{u}_1 ).
    \end{split}
\end{equation}
Therefore, the global matrices for this system are
$$\mathbf{M} = \begin{bmatrix} M_1 & \boldsymbol{0} \\\boldsymbol{0} & M_2 \end{bmatrix}, \quad \text{and}\quad \mathbf{K} = \begin{bmatrix} K^s_1 - K_1^R & K_1^L \\ K_2^L & K^s_2 - K_2^R \end{bmatrix},$$
where, in this simplified analysis based on uniform meshes, $M_1=M_2$ and $K^s_1=K^s_2$ and they are equal to the definitions
given in equation \eqref{eq:mass_stiff_mat}. The global vector of modal degrees of freedom is simply $\mathbf{U}(t)=[ \mathbf{u}_1(t) \;  \mathbf{u}_2(t)]^T$. 
Following the notation given in equation \eqref{eq:interface_mat} and figure \ref{fig:periodicDisegno}, interface matrices are defined as 
$$K^R_1 =  \phi_i^1(x_{3/2}) \phi_j^1(x_{3/2})  ,\qquad K^L_1 = \phi_i^1(x_{1/2}) \phi_j^{2}(x_{5/2}),$$
and  
$$K^R_2 =  \phi_i^2(x_{5/2}) \phi_j^2(x_{5/2})  ,\qquad K^L_2 = \phi_i^2(x_{3/2}) \phi_j^{1}(x_{3/2}).$$
%
It can be noticed that, with periodic boundary conditions, the standard CFL restriction for DG is obtained: meaning $\text{CFL}^{1}_{max}\approx 0.33$, 
$\text{CFL}^{2}_{max} \approx 0.2$, and $\text{CFL}^{3}_{max} \approx 0.14$, which approximately corresponds 
to $\text{CFL}^p_{max} \approx\frac{1}{2p+1}$.
The maximum amplification factor for periodic boundary conditions is presented in figure \ref{fig:Amplification_periodic_explicit}.

In the following analysis, a normalized $\text{CFL}\in [0,1]$ is considered by taking into account the maximum values obtained in this analysis, 
i.e. $\frac{\Delta t}{\Delta x} = \text{CFL}^{p}_{max} \text{CFL} $.
This will simplify the visualization of the stability regions in the following sections.


\subsection{Stability analysis: explicit time integration with the shifted boundary polynomial correction}\label{sec:explicitSB}

In this section, the homogeneous Dirichlet boundary condition is considered to be applied on a point that does not
coincide with the left interface of the first (boundary) element. Since an upwind flux is considered, there is no need to impose
any boundary condition of the right interface of the second element. Recall that the analysis only includes two elements to simplify
the eigenvalue spectrum visualization.

Following the notation in section \ref{sec:SBcorrection}, the analysis will consider both situations shown in figure \ref{fig:boundary}
where the boundary condition is {\it shifted} to the left interface of the first cell (i.e.\ $d<0$ and $d>0$, with $d\in[-\Delta x,\Delta x]$). 
For simplicity, the analysis is performed with $\Delta x=1$ thus the distance will be varied considering $d\in[-1,1]$.
Having either $d<0$ or $d>0$ is a choice that can be taken when meshing a domain and embedding an internal boundary: $d<0$ means that the element crossed by the boundary is not considered as an active element, while $d>0$ means that the element crossed by the boundary is kept as an active element of the discretization. Herein, we study both configurations to illustrate both possibilities. 
Due to the upwind flux, the left interface will see the flux given by the shifted boundary polynomial correction 
\begin{equation}\label{eq:SBcorrectionAnalysis}
u^\star(x_{1/2}) = 0 - ( u(x_{1/2}+d) - u(x_{1/2}) ) =  0 - \sum_{j=0}^p ( \phi^1_j(x_{1/2}+d) - \phi^1_j(x_{1/2}) )  u_{1,j},
\end{equation}
where a homogeneous Dirichlet condition, $u_D=0$, is imposed on $\bar x$,
and $u(x_{1/2}+d)$ is the value of the solution of the boundary cell on the real boundary $\bar x$.
Notice that $\phi^1_j(x_{1/2}+d)$, for $d<0$, corresponds to the evaluation of the solution outside of the boundary cell, which shows the extrapolation performed by the shifted boundary approach.

In this case, the global system is
\begin{equation}
\begin{split}
	M_1 \frac{\text{d} \mathbf{u}_1}{\text{d} t} &= K^s_1 \mathbf{u}_1 - ( K^R_1 \mathbf{u}_1 - K^{SB}_1 \mathbf{u}_1 ), \\
	M_2 \frac{\text{d} \mathbf{u}_2}{\text{d} t} &= K^s_2 \mathbf{u}_2 - ( K^R_2 \mathbf{u}_2 - K^L_2 \mathbf{u}_1 ),
\end{split}
\end{equation}
where the left interface matrix that introduces the shifted boundary polynomial correction is defined as
$$ K^{SB}_1 = -\phi^1_i(x_{1/2}) ( \phi^1_j(x_{1/2}+d) - \phi^1_j(x_{1/2}) ). $$

Therefore, the global matrices for this system are
\begin{equation}\label{eq:SBexplicitMatrix}
\mathbf{M} = \begin{bmatrix} M_1 & \boldsymbol{0} \\\boldsymbol{0} & M_2 \end{bmatrix}, \quad \text{and}\quad \mathbf{K} = \begin{bmatrix} K^s_1 - K_1^R + K_1^{SB} & \boldsymbol{0} \\ K_2^L & K^s_2 - K_2^R \end{bmatrix},
\end{equation}
where the other local matrices are defined as in section \ref{sec:periodic}.

\begin{figure}[h!]
\centering
\includegraphics[width=0.95\textwidth]{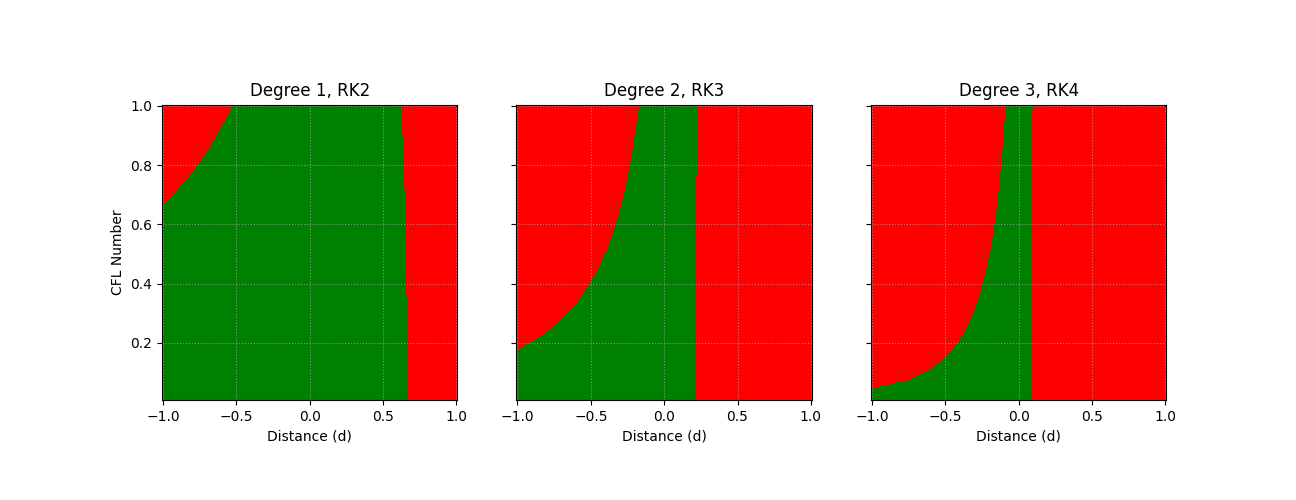}
\caption{Stability region: shifted boundary polynomial correction of the homogeneous Dirichlet condition for the explicit discontinuous Galerkin method. Green areas are stable, while red areas are unstable.}
\label{fig:DG_SBM2_explicit}
\end{figure}

The stability region for explicit discontinuous Galerkin schemes with shifted boundary polynomial correction is shown in figure \ref{fig:DG_SBM2_explicit}.
As in the previous section, the stability region is given by the amplification factor computed as the spectral radius of the eigenvalues of the discretized operators, which are amplified with the corresponding high order time integration method \eqref{eq:amplificationRK}. 
It should be noticed that the stability region is not symmetric with respect to the value of $d\in[-1,1]$. 
When $d<0$, there are stable areas for all values of $d$ with a constraint on the CFL number that becomes stricter as $d\rightarrow-1$ and $p$ increases, 
while for $d>0$ the stability region of the method is not constrained by the CFL number but the method becomes 
unstable after a certain threshold $d>d_{max}$, which depends on the polynomial degree used to build the scheme.
In particular, when $d=-1$, the maximum admissible time step is dictated by a $\text{CFL}_{max}\approx 0.6$ for $\mathbb{P}^1$,
$\text{CFL}_{max}\approx 0.2$ for $\mathbb{P}^2$, and $\text{CFL}_{max}\approx 0.05$ for $\mathbb{P}^3$.
Moreover, when $d>0$, $d_{max}\approx 0.6$ for $\mathbb{P}^1$, $d_{max}\approx 0.2$ for $\mathbb{P}^2$, and $d_{max}\approx 0.1$ for $\mathbb{P}^3$. 

\begin{figure}[h!]
\centering
\subfigure[$\mathbb{P}^1$, $d<0$]{\includegraphics[width=0.45\textwidth]{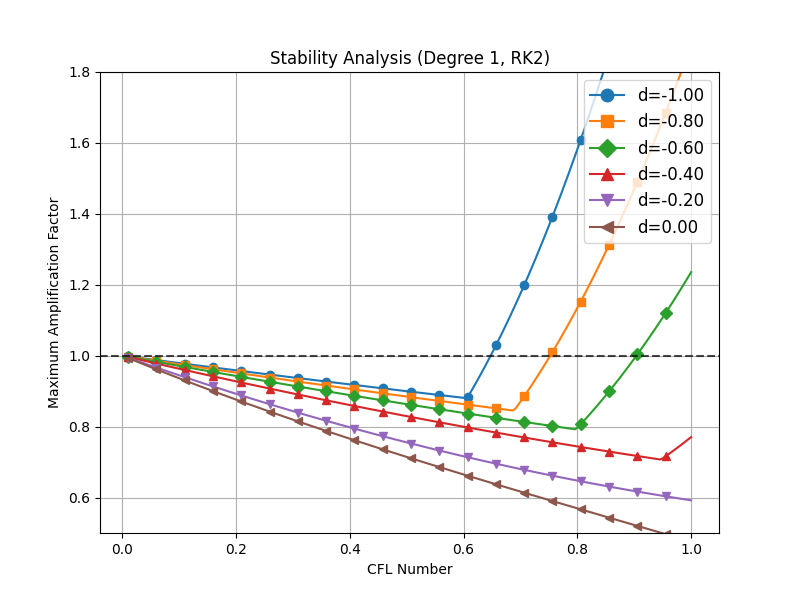}}
\subfigure[$\mathbb{P}^1$, $d>0$]{\includegraphics[width=0.45\textwidth]{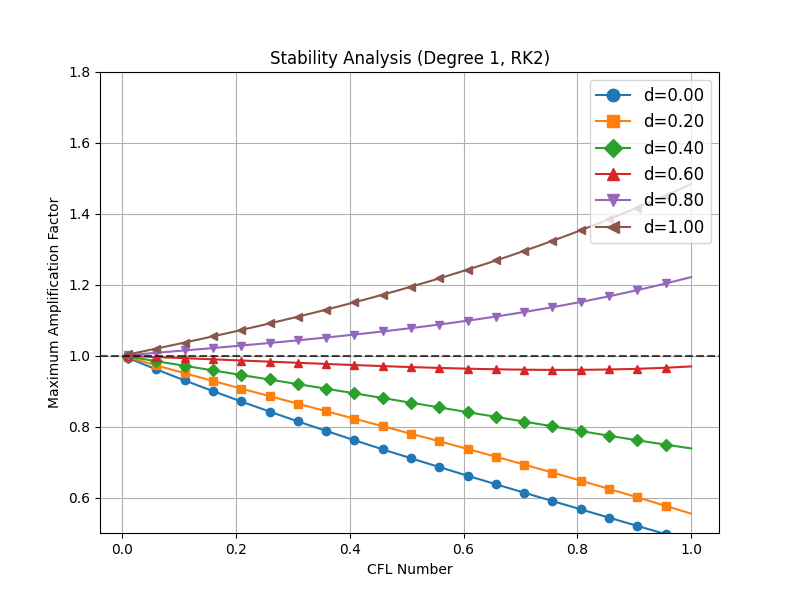}}\quad
\subfigure[$\mathbb{P}^2$, $d<0$]{\includegraphics[width=0.45\textwidth]{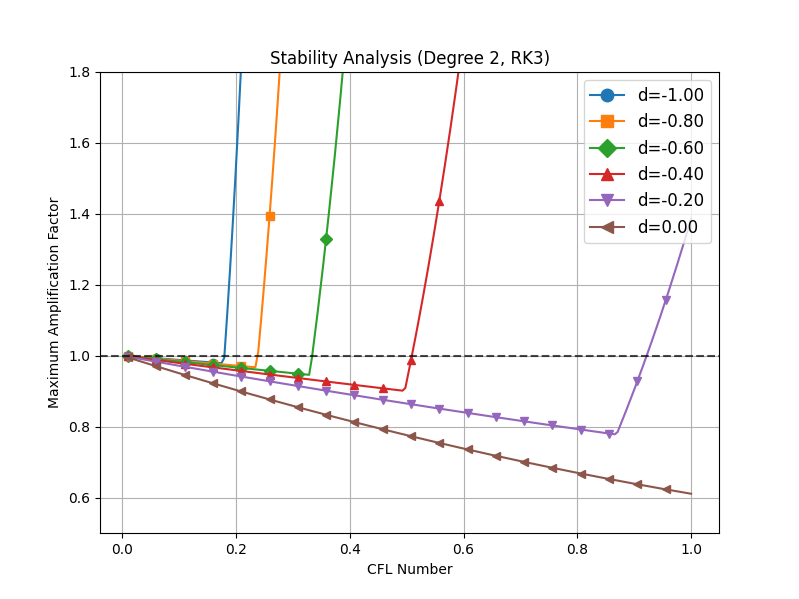}}
\subfigure[$\mathbb{P}^2$, $d>0$]{\includegraphics[width=0.45\textwidth]{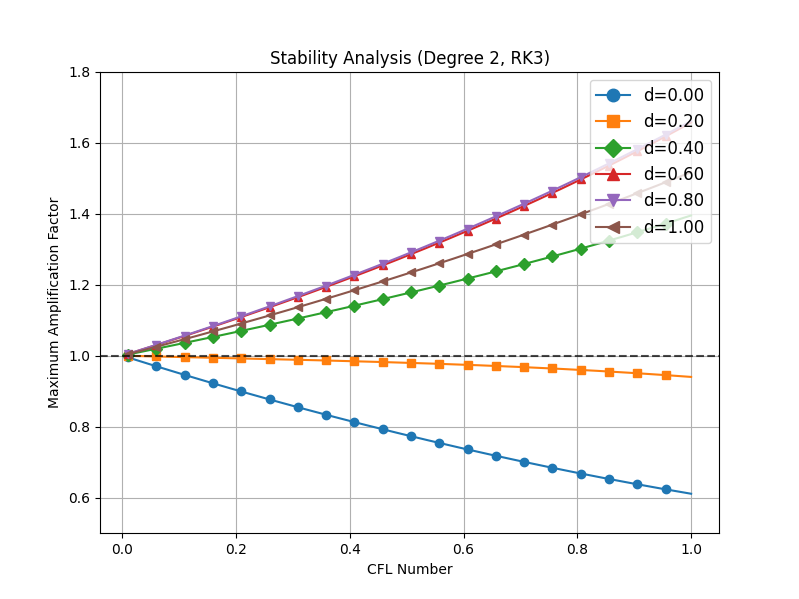}}\quad
\subfigure[$\mathbb{P}^3$, $d<0$]{\includegraphics[width=0.45\textwidth]{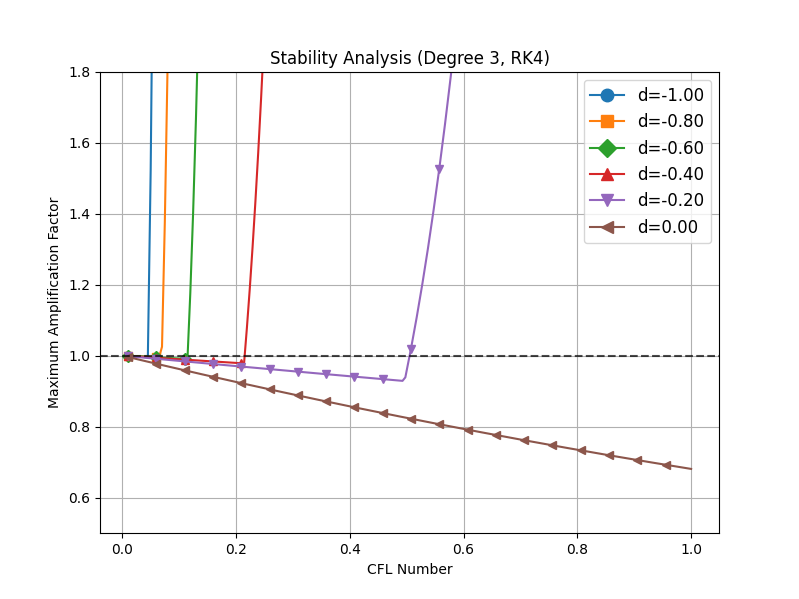}}
\subfigure[$\mathbb{P}^3$, $d>0$]{\includegraphics[width=0.45\textwidth]{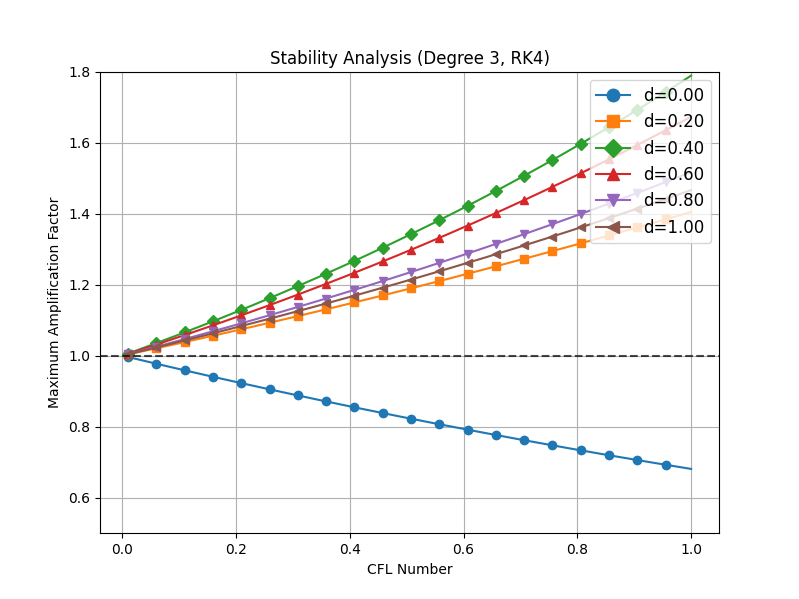}}
\caption{Amplification factor: shifted boundary polynomial correction of the homogeneous Dirichlet condition for the explicit discontinuous Galerkin method.}
\label{fig:Amplification_SBM2_explicit}
\end{figure}

The maximum amplification factor is also presented in figure \ref{fig:Amplification_SBM2_explicit} for different values of $d<0$ and $d>0$.
For the former case, it can be observed that by increasing the polynomial order, the amplification factor move faster to the unstable region and
a lower CFL value is needed to preserve stability. Instead, for the latter, it is possible to notice how the scheme becomes
unconditionally unstable for $d\rightarrow 1$, while staying stable under the standard CFL condition for $d<d_{max}$. 

It is also interesting to look at the characteristic polynomial and analytical eigenvalues of the semi-discrete system. 
Indeed, thanks to the characteristic polynomial, it is possible to prove that only the eigenvalues of the boundary cell are affected by the shifted boundary polynomial correction, i.e.\ they depend on the distance $d$,
while the eigenvalues of the internal cells (no matter on many there are) do not depend on $d$.
For a two cells system with a $\mathbb{P}^1$ method, the characteristic polynomial reads
$$ \mathcal{P}_{\mathbb{P}^1}(\lambda)=\frac{\left(\lambda^{2} + 4 \lambda + 6\right) \left(\lambda^{2} - 6 \lambda d + 4 \lambda + 6\right)}{12^2}. $$
By increasing the number of cells, for example to three cells, the characteristic polynomial becomes
$$ \mathcal{P}_{\mathbb{P}^1}(\lambda)= \frac{\left(\lambda^{2} + 4 \lambda + 6\right)^{2} \left(\lambda^{2} - 6 \lambda d + 4 \lambda + 6\right)}{12^3},$$
and so on for more cells.

Similarly for $\mathbb{P}^2$ the characteristic polynomial for a two cells system is
$$ \mathcal{P}_{\mathbb{P}^2}(\lambda) = \frac{\left(\lambda^{3} + 9 \lambda^{2} + 36 \lambda + 60\right) \left(\lambda^{3} + 30 \lambda^{2} d^{2} - 36 \lambda^{2} d + 9 \lambda^{2} - 60 \lambda d + 36 \lambda + 60\right)}{135^2}, $$
while for a three cells system it becomes
$$ \mathcal{P}_{\mathbb{P}^2}(\lambda) = \frac{\left(\lambda^{3} + 9 \lambda^{2} + 36 \lambda + 60\right)^{2} \left(\lambda^{3} + 30 \lambda^{2} d^{2} - 36 \lambda^{2} d + 9 \lambda^{2} - 60 \lambda d + 36 \lambda + 60\right)}{135^3}. $$
For $\mathbb{P}^3$, the characteristic polynomial becomes extremely cumbersome to be reported here.

For $\mathbb{P}^1$, analytical eigenvalues can be easily computed by hand, and read as
\begin{equation}\label{eq:eigen}
\begin{split}
	\lambda_1 &= -2 - \sqrt{2} \im \\
	\lambda_2 &= -2 + \sqrt{2} \im \\
	\lambda_3 &= 3 d - \sqrt{9 d^{2} - 12 d - 2} - 2 \\
	\lambda_4 &= 3 d + \sqrt{9 d^{2} - 12 d - 2} - 2 
\end{split}
\end{equation}
where $\lambda_1$ and $\lambda_2$ are the eigenvalues of the cell far from the boundary, 
while $\lambda_3$ and $\lambda_4$ are the ones of the cell affected by the shifted boundary polynomial correction.
As expected, the eigenvalues depend on the distance $d$ of the immersed boundary from the left edge of the cell,
which may bring either a positive or negative real part. 
In particular, the real part of $\lambda_3$ and $\lambda_4$ is negative if $d<2/3\approx 0.66$ which is the distance threshold after
which the $\mathbb{P}^1$ scheme becomes unconditionally unstable.
For $\mathbb{P}^2$  and $\mathbb{P}^3$, the symbolic algebra for computing analytically symbolic eigenvalues becomes intractable,
thus only the numerical eigenvalues are trusted. 

\subsection{Stability analysis: implicit time integration with the shifted boundary polynomial correction}\label{sec:implicitSB}

Due to the strong constraints imposed by explicit time integration, in this section, the possibility of performing implicit integration for the discontinuous Galerkin method with the shifted boundary polynomial correction is explored.
Since the analysis mainly focuses on how the embedded boundary treatment affects the stability region of a high order DG space discretization, only implicit Euler time integration is considered to study if it is possible to relax such constraints. The extension of this method to higher order implicit schemes \cite{offner2025analysis} will be developed in a different work, which will focus more on the accuracy for time-dependent problems with embedded configurations. 

The accuracy in the numerical experiments presented in section \ref{sec:tests} will not be compromised simply because
a stationary manufactured solution will be considered.

If we consider the implicit Euler time integration of the discontinuous Galerkin method with the shifted boundary polynomial correction,
the global system reads

\begin{equation}
\begin{split}
M_1\frac{\mathbf{u}^{n+1}_1 - \mathbf{u}^{n}_1}{\Delta t} &= K^s_1 \mathbf{u}^{n+1}_1 - ( K^R_1 \mathbf{u}^{n+1}_1 - K^{SB}_1 \mathbf{u}^{n+1}_1),\\
M_2\frac{\mathbf{u}^{n+1}_2 - \mathbf{u}^{n}_2}{\Delta t} &= K^s_2 \mathbf{u}^{n+1}_2 - ( K^R_2 \mathbf{u}^{n+1}_2 - K^L_2 \mathbf{u}^{n+1}_1 ),
\end{split}
\end{equation}

In this case, the system reduces to find the eigenvalues of the following amplification matrix:
\begin{equation}
\begin{bmatrix} \mathbf{u}^{n+1}_1 \\ \mathbf{u}^{n+1}_2 \end{bmatrix} = 
\left[I - \Delta t \mathbf{M}^{-1} \mathbf{K}\right]^{-1}
\begin{bmatrix} \mathbf{u}^{n}_1 \\ \mathbf{u}^{n}_2 \end{bmatrix}
\end{equation}
where the global matrices are those defined in equation \eqref{eq:SBexplicitMatrix}.

\begin{figure}[h!]
\centering
\includegraphics[width=0.95\textwidth]{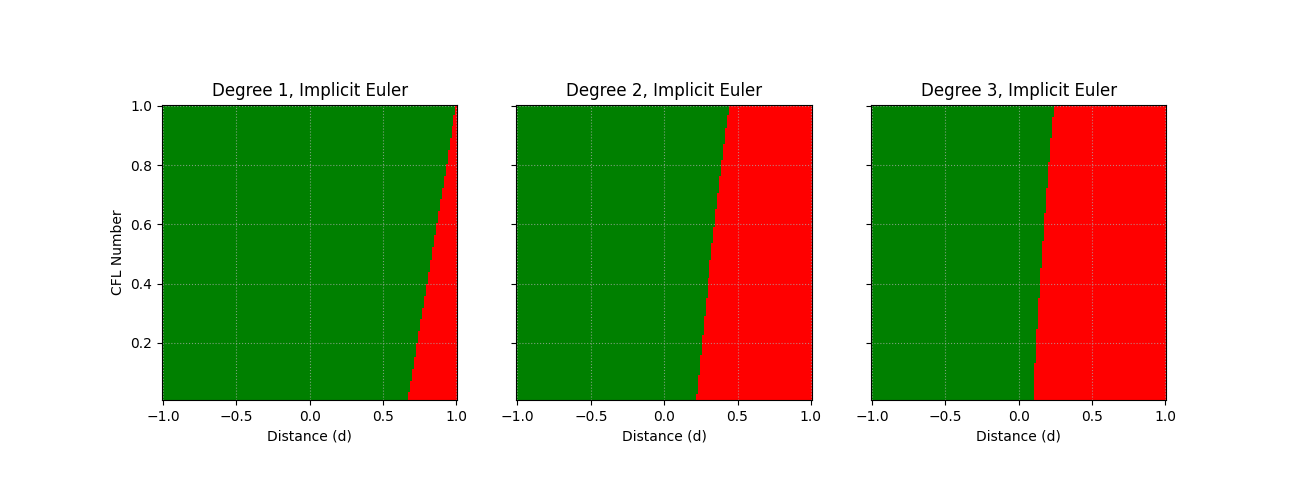}
\caption{Stability region: shifted boundary polynomial correction of the homogeneous Dirichlet condition for the implicit discontinuous Galerkin method with $\text{CFL}\in[0,1]$. Green areas are stable, while red areas are unstable.}
\label{fig:DG_SBM2_implicit1}
\end{figure}

In figure \ref{fig:DG_SBM2_implicit1}, the stability region for the implicit time discretization is visualized when considering $\text{CFL}\in[0,1]$.
It should be noticed that also for the implicit case, the stability region is not symmetric with respect to $d\in[-1,1]$. 
When $d<0$, there are stable areas for all values of $d$ with no constraint on the CFL number, 
which is an important improvement with respect to the explicit case.
When considering the range $\text{CFL}\in[0,1]$ and $d>0$, the stability region of the method is not constrained by the CFL number 
but the method becomes unconditionally unstable after a certain threshold $d>d_{max}$, as also observed in the explicit case.
In this case, $d_{max}\approx 0.9$ for $\mathbb{P}^1$, $d_{max}\approx 0.4$ for $\mathbb{P}^2$, and $d_{max}\approx 0.2$ for $\mathbb{P}^3$.
Recall that, although in this analysis $d$ can be either positive or negative, for practical applications the choice of the position of the embedded geometry
(see figure \ref{fig:boundary}) can be taken always negative, resp.\ positive, by simply eliminating, resp.\ keeping, the elements cut by the boundary.  
Therefore, if one is able to always eliminating the cut elements, the distance will always be negative, and the scheme unconditionally stable.

\begin{figure}[h!]
\centering
\subfigure[$\mathbb{P}^1$, $d<0$]{\includegraphics[width=0.45\textwidth]{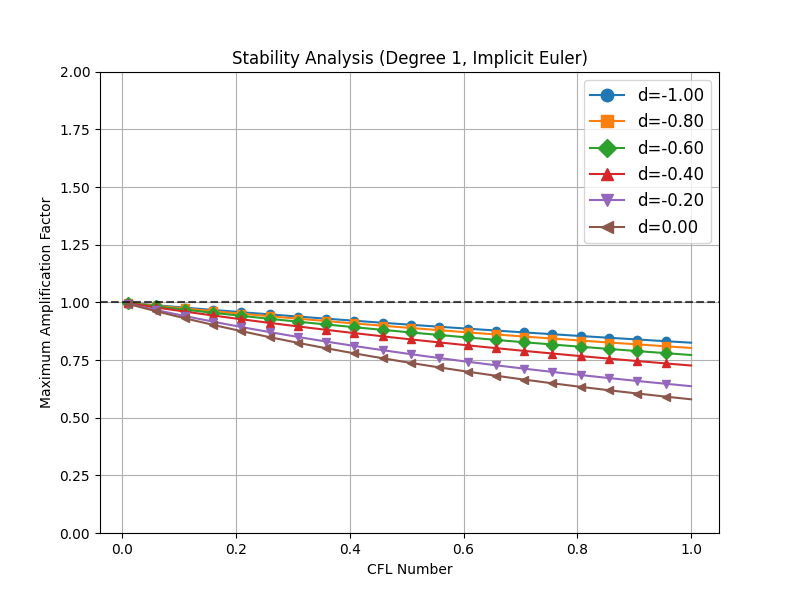}}
\subfigure[$\mathbb{P}^1$, $d>0$]{\includegraphics[width=0.45\textwidth]{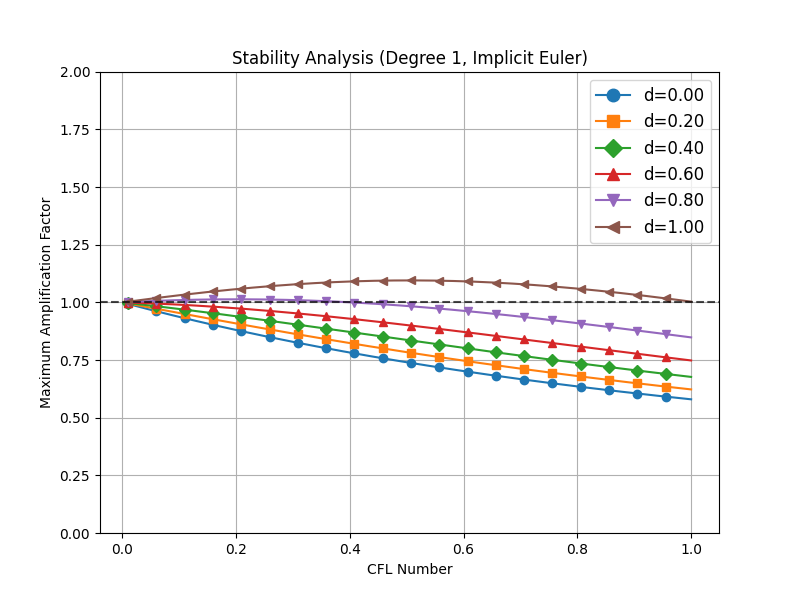}}\quad
\subfigure[$\mathbb{P}^2$, $d<0$]{\includegraphics[width=0.45\textwidth]{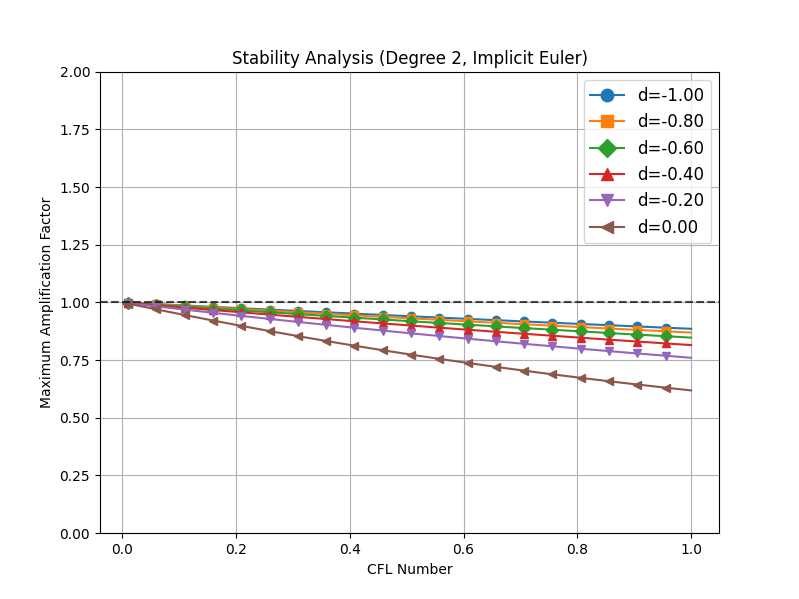}}
\subfigure[$\mathbb{P}^2$, $d>0$]{\includegraphics[width=0.45\textwidth]{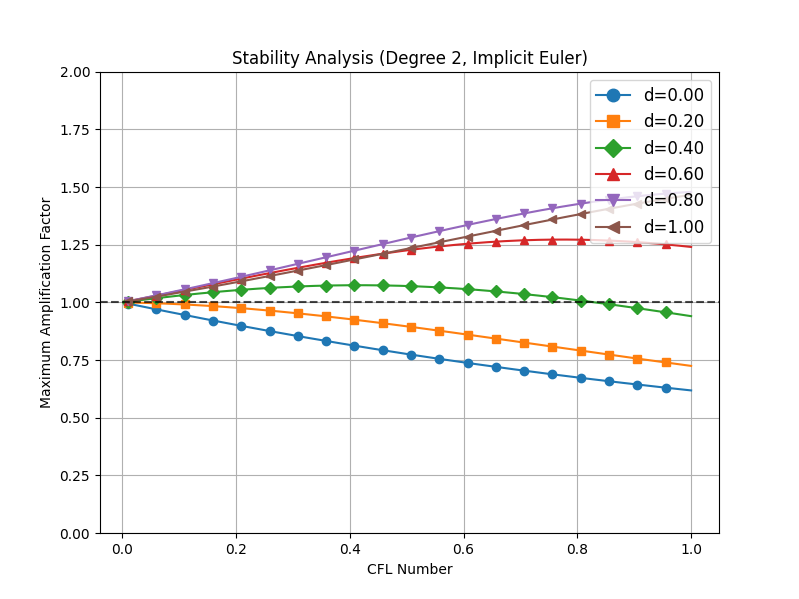}}\quad
\subfigure[$\mathbb{P}^3$, $d<0$]{\includegraphics[width=0.45\textwidth]{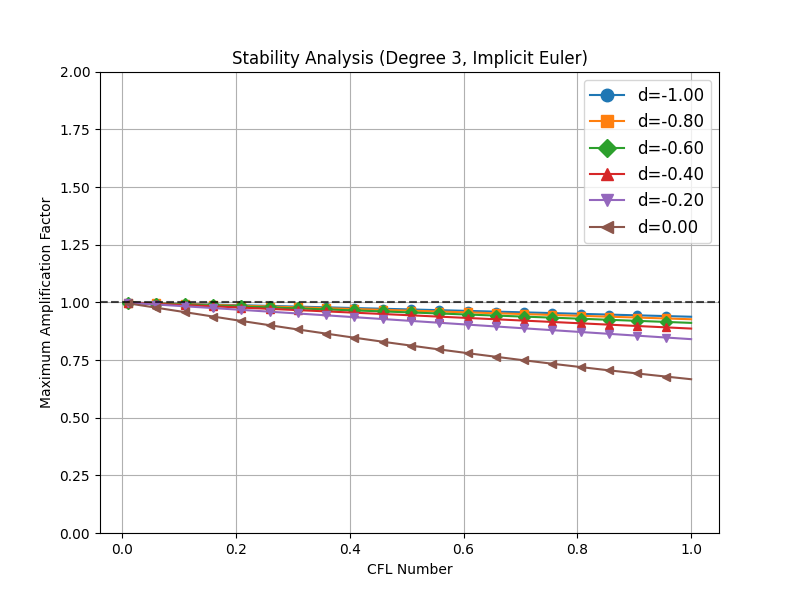}}
\subfigure[$\mathbb{P}^3$, $d>0$]{\includegraphics[width=0.45\textwidth]{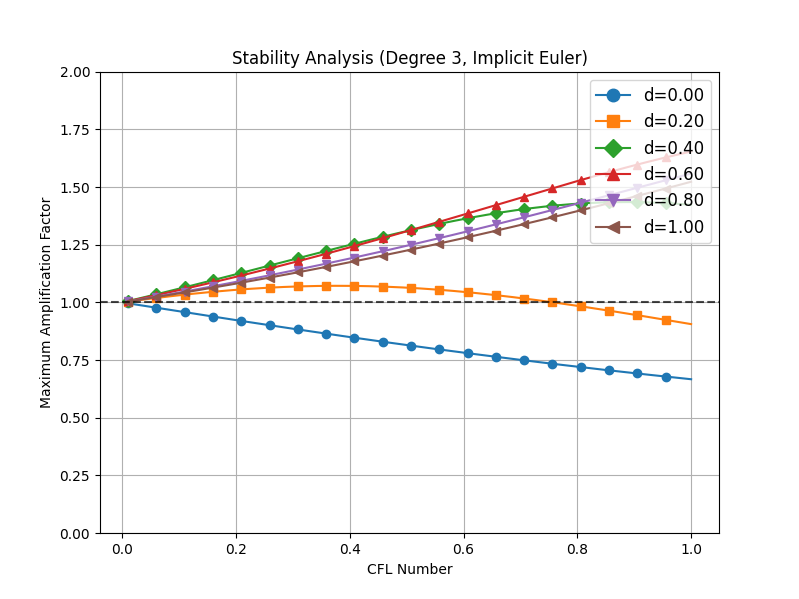}}
\caption{Amplification factor: shifted boundary polynomial correction of the homogeneous Dirichlet condition for the implicit discontinuous Galerkin method with $\text{CFL}\in[0,1]$.}
\label{fig:Amplification_SBM2_implicit1}
\end{figure}

More details about this behavior are given by the maximum amplification factor shown in figure \ref{fig:Amplification_SBM2_implicit1}.

\begin{figure}[h!]
\centering
\includegraphics[width=0.95\textwidth]{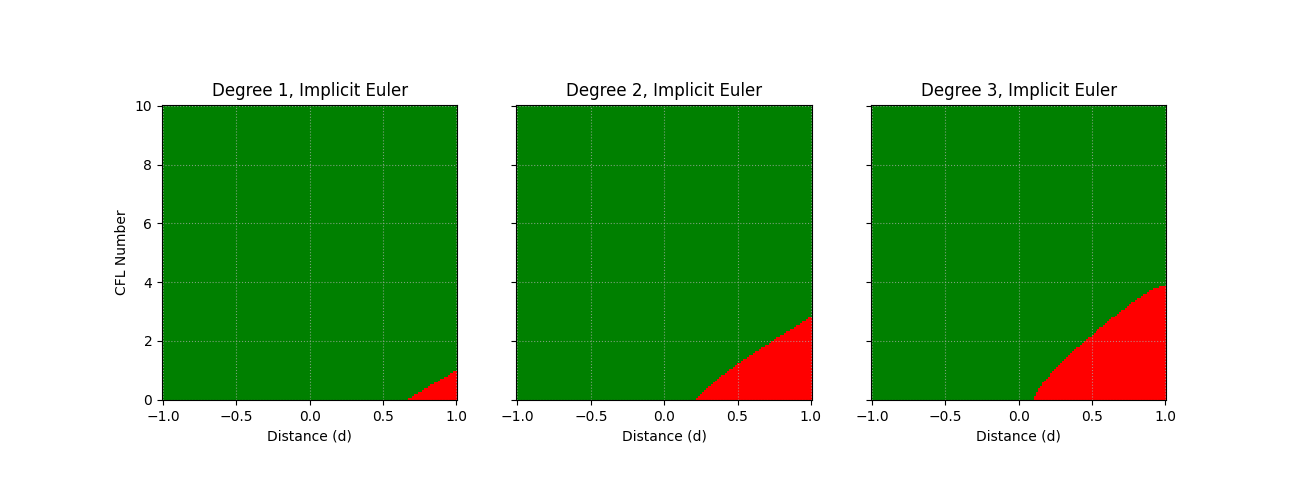}
\caption{Stability region: shifted boundary polynomial correction of the homogeneous Dirichlet condition for the implicit discontinuous Galerkin method with $\text{CFL}\in[0,10]$. Green areas are stable, while red areas are unstable.}
\label{fig:DG_SBM2_implicit10}
\end{figure}

\begin{figure}[h!]
\centering
\subfigure[$\mathbb{P}^1$, $d<0$]{\includegraphics[width=0.45\textwidth]{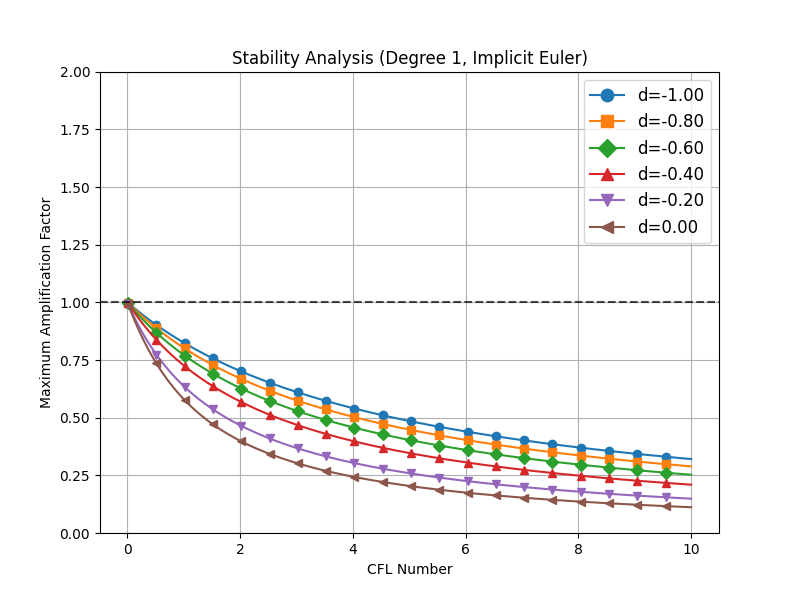}}
\subfigure[$\mathbb{P}^1$, $d>0$]{\includegraphics[width=0.45\textwidth]{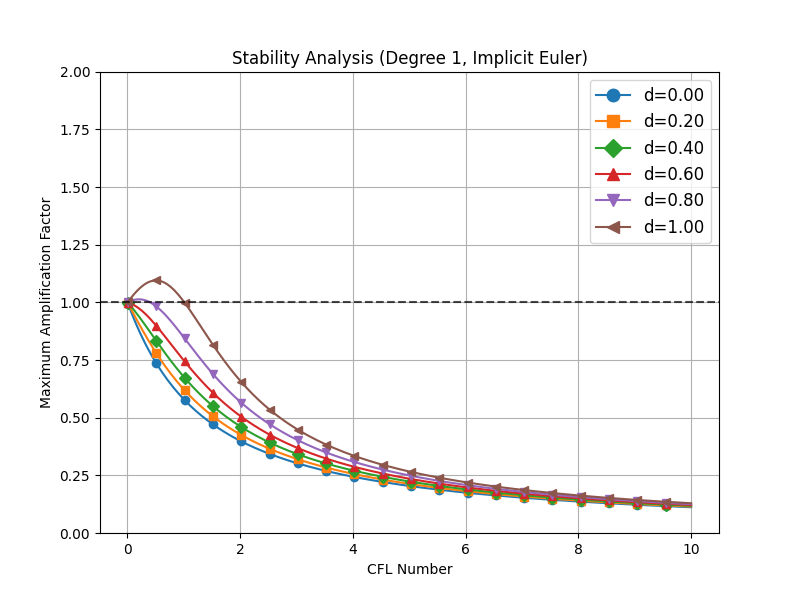}}\quad
\subfigure[$\mathbb{P}^2$, $d<0$]{\includegraphics[width=0.45\textwidth]{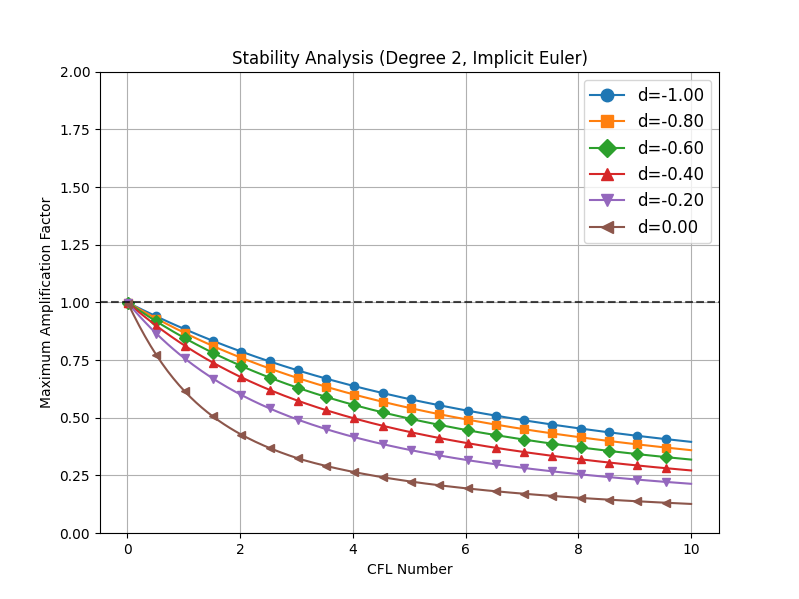}}
\subfigure[$\mathbb{P}^2$, $d>0$]{\includegraphics[width=0.45\textwidth]{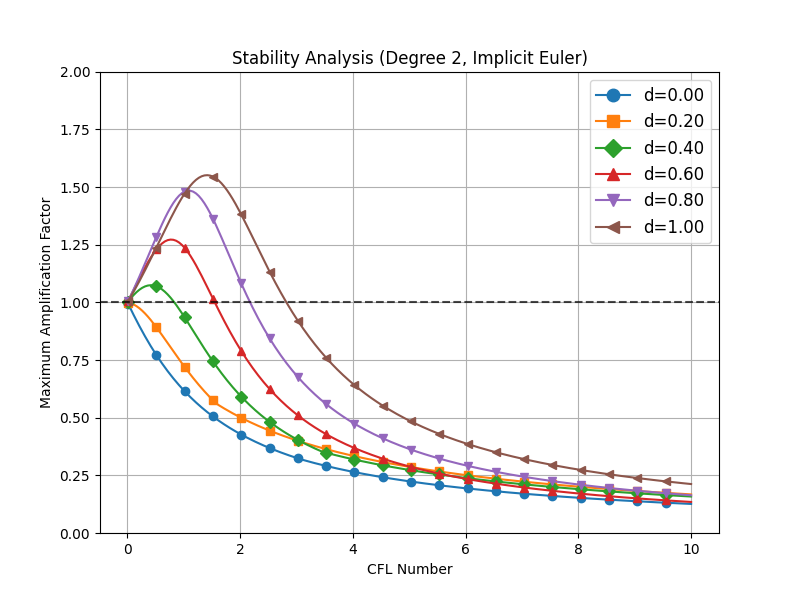}}\quad
\subfigure[$\mathbb{P}^3$, $d<0$]{\includegraphics[width=0.45\textwidth]{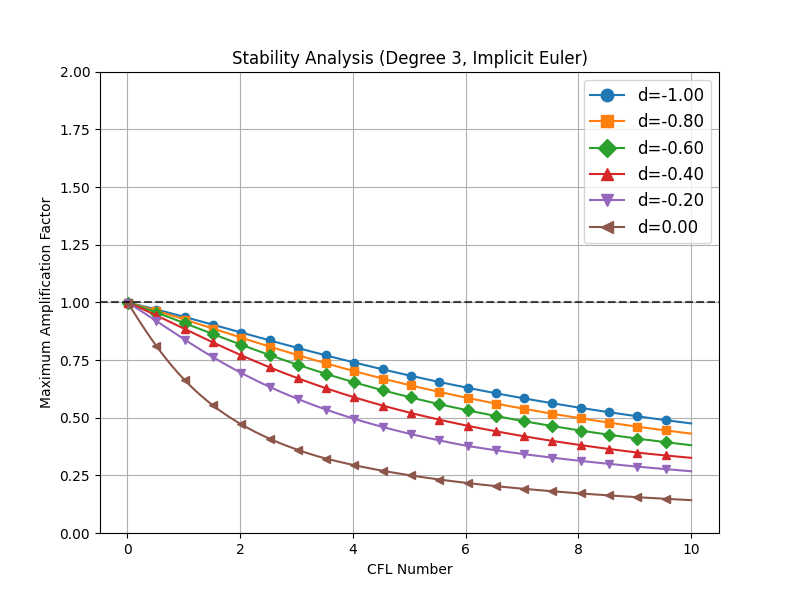}}
\subfigure[$\mathbb{P}^3$, $d>0$]{\includegraphics[width=0.45\textwidth]{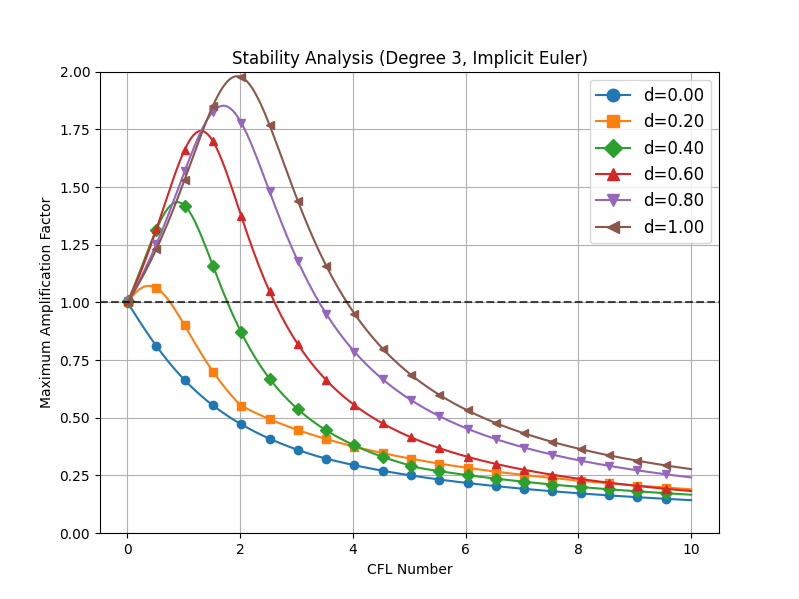}}
\caption{Amplification factor: shifted boundary polynomial correction of the homogeneous Dirichlet condition for the implicit discontinuous Galerkin method with $\text{CFL}\in[0,10]$.}
\label{fig:Amplification_SBM2_implicit10}
\end{figure}

It is also interesting to observe that, when considering $\text{CFL}\in[0,10]$, the stability region becomes much larger as also shown in figure \ref{fig:DG_SBM2_implicit10}.
In particular, in this case, there is a certain threshold of the CFL after which the scheme is unconditionally stable for all values of $d\in[-1,1]$.
In particular, when $d=1$, unconditional stability is given for a $\text{CFL}>1$ for $\mathbb{P}^1$,
$\text{CFL}>3$ for $\mathbb{P}^2$, and $\text{CFL}>4$ for $\mathbb{P}^3$.
The maximum amplification for $\text{CFL}\in[0,10]$ is given in figure \ref{fig:Amplification_SBM2_implicit10}, 
where it is clear how it moves from the unconditional instability to the unconditional stability.

\section{Numerical experiments}\label{sec:tests}

For the numerical experiments, the linear advection equation with source term is set up to preserve a stationary manufactured solution
given by the following exact solution
$$ u_{ex}(x) = 0.1 \sin(\pi x) ,$$
which is also taken as initial condition.
The source term is computed as
$$\partial_t u_{ex}(x) + \partial_x u_{ex}(x) = s(x), \quad\text{that gives}\quad s(x) = 0.1 \pi \cos(\pi x) .$$
The computational domain is $\Omega=[0,2]$ with a Dirichlet boundary condition imposed at the left boundary,
and no conditions are imposed at the right boundary due to the pure upwind flux.
A set of uniform meshes with $N_e=20,40,80,160,320$ elements is considered. 
As in the analysis the polynomials considered are $\mathbb{P}^p$ with $p=1,2,3$.
The shifted boundary polynomial correction is imposed considering different distances $d$ of the boundary from the left interface of the first cell.
When $d=0$, the boundary coincides with the left interface of the first cell, and therefore
the correction is not applied. This is the case of standard fitted boundary conditions.
The distances considered in the numerical experiments vary as $d\in[-\Delta x,\Delta x]$ 
to explore all the behaviors of the method studied in the analysis.
For simplicity, in convergence tables, the distance is reported in non-dimensional form, i.e.\ $d/\Delta x$.
Stability results obtained through the analysis are perfectly corroborated by the numerical results for all configurations of $d$,
and for explicit and implicit time integration.

In table \ref{tab:LAE_standardBC}, the results obtained with explicit time integration and fitted boundary conditions are presented.
The fitted results, which validate the analysis in section \ref{sec:periodic}, are only shown here as a baseline for embedded configurations.

\begin{table}
        \caption{Linear advection equation: $L_2$ errors $\|u-u_{ex}\|_2$ and estimated order of accuracy (EOA) with fitted boundary condition and explicit time integration 
		with distance $d=0$, CFL$=1$, and polynomials $\mathbb{P}^p$.}\label{tab:LAE_standardBC}
        \scriptsize
        \centering
        \begin{tabular}{c|cc|cc|cc} \hline\hline
                        &\multicolumn{2}{c|}{$\mathbb{P}^1$} &\multicolumn{2}{c|}{$\mathbb{P}^2$} &\multicolumn{2}{c}{$\mathbb{P}^3$}  \\[0.5mm]
                        \cline{2-7}
                        $N_e$ & $L_2$ error & EOA & $L_2$  error       & EOA & $L_2$ error       &EOA  \\ \hline 
						20    &  4.75E-04   & --   &  1.16E-05  &  --    &  2.19E-07  &   --     \\
                        40    &  1.19E-04   & 2.00 &  1.44E-06  &  3.00  &  1.37E-08  &  4.00      \\
                        80    &  2.97E-05   & 2.00 &  1.81E-07  &  3.00  &  8.52E-10  &  4.00      \\
                        160   &  7.42E-06   & 2.00 &  2.26E-08  &  3.00  &  5.00E-11  &  4.09      \\
                        320   &  1.85E-06   & 2.00 &  2.82E-09  &  3.00  &  4.43E-12  &  3.49      \\ \hline
        \end{tabular}
\end{table}

In table \ref{tab:LAE_SBexplicit}, the numerical results obtained for different embedded configurations are shown in the case
of explicit time integration. As can be noticed, when taking the constrained CFL for $d<0$ or the maximum distance threshold $0<d<d_{max}$,
all convergence slopes are recovered. For all the couples of parameters $(\text{CFL},d)$ for which the analysis predicted instabilities,
the numerical solutions are indeed unstable, and are therefore not shown here.  

\begin{table}
        \caption{Linear advection equation: $L_2$ errors $\|u-u_{ex}\|_2$ and estimated order of accuracy (EOA) with shifted boundary polynomial correction and explicit time integration 
		by varying the distance $d$ for polynomials $\mathbb{P}^p$.}\label{tab:LAE_SBexplicit}
        \scriptsize
        \centering
        \begin{tabular}{c|cc|cc|cc} \hline\hline
                        &\multicolumn{2}{c|}{$d=-1.0$, CFL$=0.6$} &\multicolumn{2}{c|}{$d=-0.5$, CFL$=1.0$} &\multicolumn{2}{c}{$d=+0.6$, CFL$=1.0$}  \\[0.5mm]
                        \cline{2-7}
                        $N_e$, $\mathbb{P}^1$ & $L_2$ error        & EOA & $L_2$  error       & EOA & $L_2$ error       &EOA  \\ \hline 
                        20    &  5.98E-04   &  --  &  6.57E-04  &  --   &  6.87E-04  &  --      \\
                        40    &  1.27E-04   & 2.23 &  1.32E-04  & 2.31  &  1.34E-04  &   2.35   \\
                        80    &  3.02E-05   & 2.07 &  3.05E-05  & 2.10  &  3.07E-05  &   2.12   \\
                        160   &  7.45E-06   & 2.01 &  7.47E-06  & 2.03  &  7.48E-06  &   2.03   \\
                        320   &  1.86E-06   & 2.00 &  1.86E-06  & 2.00  &  1.86E-06  &   2.00   \\ \hline
        \end{tabular}
        \begin{tabular}{c|cc|cc|cc} 
                        &\multicolumn{2}{c|}{$d=-1.0$, CFL$=0.2$} &\multicolumn{2}{c|}{$d=-0.5$, CFL$=0.4$} &\multicolumn{2}{c}{$d=+0.2$, CFL$=1.0$}  \\[0.5mm]
                        \cline{2-7}
                        $N_e$, $\mathbb{P}^2$ & $L_2$ error        & EOA & $L_2$  error       & EOA & $L_2$ error       &EOA  \\ \hline 
                        20    &  2.68E-03   & --    &  7.43E-04    &    --  &  8.45E-05  &   --    \\
                        40    &  3.37E-04   &  3.00 &  9.34E-05    &  3.00  &  1.07E-05  &  2.98     \\
                        80    &  4.22E-05   &  3.00 &  1.17E-05    &  3.00  &  1.34E-06  &  3.00     \\
                        160   &  5.28E-06   &  3.00 &  1.46E-06    &  3.00  &  1.67E-07  &  3.00     \\
                        320   &  6.60E-07   &  3.00 &  1.83E-07    &  3.00  &  2.09E-08  &  3.00     \\   \hline
        \end{tabular}
        \centering
        \begin{tabular}{c|cc|cc|cc} 
                        &\multicolumn{2}{c|}{$d=-1.0$, CFL$=0.05$} &\multicolumn{2}{c|}{$d=-0.5$, CFL$=0.15$} &\multicolumn{2}{c}{$d=+0.1$, CFL$=1.0$}  \\[0.5mm]
                        \cline{2-7}
                        $N_e$, $\mathbb{P}^3$ & $L_2$ error     & EOA & $L_2$  error       & EOA & $L_2$ error       &EOA  \\ \hline 
                        20    &  2.59E-05   & --   &  6.62E-06  & --   &  3.46E-07  &  --   \\
                        40    &  8.12E-07   & 4.99 &  2.08E-07  & 5.00 &  1.61E-08  &  4.42     \\
                        80    &  2.54E-08   & 5.00 &  6.53E-09  & 5.00 &  8.92E-10  &  4.17     \\
                        160   &  7.65E-10   & 5.04 &  2.01E-10  & 5.02 &  5.06E-11  &  4.14     \\
                        320   &  1.07E-11   & 6.16 &  4.99E-12  & 5.33 &  4.43E-12  &  3.51     \\   \hline
                        \hline
        \end{tabular}
\end{table}

Table \ref{tab:LAE_SBimplicit1} presents the results obtained for embedded configurations with the implicit time integration, 
and considering for all cases $\text{CFL}=1$. For $d<0$ the scheme is always stable.
For $d>0$ the scheme is stable only when $0<d<d_{max}$ for the distance threshold described in section \ref{sec:implicitSB},
and it is indeed unstable for $d>d_{max}$. 
Owing to the high computational cost of the current, unoptimized implementation of the implicit scheme, the finest mesh is omitted from Table \ref{tab:LAE_SBimplicit1} and will be presented in Table \ref{tab:LAE_SBimplicit100}.

\begin{table}
        \caption{Linear advection equation: $L_2$ errors $\|u-u_{ex}\|_2$ and estimated order of accuracy (EOA) with shifted boundary polynomial correction and implicit time integration 
		by varying the distance $d$ for polynomials $\mathbb{P}^p$. The CFL is always taken equal to 1.}\label{tab:LAE_SBimplicit1}
        \scriptsize
        \centering
        \begin{tabular}{c|cc|cc|cc} \hline\hline
                        &\multicolumn{2}{c|}{$d=-1.0$, CFL$=1.0$} &\multicolumn{2}{c|}{$d=-0.5$, CFL$=1.0$} &\multicolumn{2}{c}{$d=+0.9$, CFL$=1.0$}  \\[0.5mm]
                        \cline{2-7}
                        $N_e$, $\mathbb{P}^1$ & $L_2$ error        & EOA & $L_2$  error       & EOA & $L_2$ error       &EOA  \\ \hline 
                        20    &   5.98E-04	&  --  & 6.57E-04 &  --  & 6.55E-04 &  --   \\
                        40    &   1.27E-04	& 2.23 & 1.32E-04 & 2.31 & 1.32E-04 & 2.31  \\
                        80    &   3.02E-05	& 2.07 & 3.05E-05 & 2.10 & 3.05E-05 & 2.10  \\
                        160   &   7.45E-06	& 2.01 & 7.47E-06 & 2.03 & 7.47E-06 & 2.02  \\ \hline
        \end{tabular}
        \begin{tabular}{c|cc|cc|cc} 
                        &\multicolumn{2}{c|}{$d=-1$, CFL$=1.0$} &\multicolumn{2}{c|}{$d=-0.5$, CFL$=1.0$} &\multicolumn{2}{c}{$d=+0.4$, CFL$=1.0$}  \\[0.5mm]
                        \cline{2-7}
                        $N_e$, $\mathbb{P}^2$ & $L_2$ error        & EOA & $L_2$  error       & EOA & $L_2$ error       &EOA  \\ \hline 
                        20    & 2.68E-03 &  --  & 7.43E-04 &  --  & 9.86E-05 &  --   \\
                        40    & 3.37E-04 & 3.00 & 9.34E-05 & 3.00 & 1.25E-05 & 2.98   \\
                        80    & 4.22E-05 & 3.00 & 1.17E-05 & 3.00 & 1.56E-06 & 3.00   \\
                        160   & 5.28E-06 & 3.00 & 1.46E-06 & 3.00 & 1.95E-07 & 3.00   \\ \hline
        \end{tabular}
        \centering
        \begin{tabular}{c|cc|cc|cc} 
                        &\multicolumn{2}{c|}{$d=-1$, CFL$=1.0$} &\multicolumn{2}{c|}{$d=-0.5$, CFL$=1.0$} &\multicolumn{2}{c}{$d=+0.2$, CFL$=1.0$}  \\[0.5mm]
                        \cline{2-7}
                        $N_e$, $\mathbb{P}^3$ & $L_2$ error     & EOA & $L_2$  error       & EOA & $L_2$ error       &EOA  \\ \hline 
                        20    & 2.59E-05  &  --  & 6.62E-06 &  --  & 4.19E-07 &  --   \\
                        40    & 8.12E-07  & 4.99 & 2.08E-07 & 5.00 & 1.77E-08 &  4.56 \\
                        80    & 2.54E-08  & 5.00 & 6.53E-09 & 5.00 & 9.21E-10 &  4.26 \\
                        160   & 7.65E-10  & 5.04 & 2.01E-10 & 5.02 & 5.10E-11 &  4.17 \\ \hline
                        \hline
        \end{tabular}
\end{table}

On the contrary, when taking $\text{CFL}=100$, the scheme is unconditionally stable for all embedded configurations as shown in table \ref{tab:LAE_SBimplicit100}.
This result completes the validation of the analysis performed in section \ref{sec:explicitSB} and \ref{sec:implicitSB}. 

\begin{table}
        \caption{Linear advection equation: $L_2$ errors $\|u-u_{ex}\|_2$ and estimated order of accuracy (EOA) with shifted boundary polynomial correction and implicit time integration 
		by varying the distance $d$ for polynomials $\mathbb{P}^p$.The CFL is always taken equal to 100.}\label{tab:LAE_SBimplicit100}
        \scriptsize
        \centering
        \begin{tabular}{c|cc|cc|cc} \hline\hline
                        &\multicolumn{2}{c|}{$d=-1.0$, CFL$=100$} &\multicolumn{2}{c|}{$d=-0.5$, CFL$=100$} &\multicolumn{2}{c}{$d=+0.8$, CFL$=100$}  \\[0.5mm]
                        \cline{2-7}
                        $N_e$, $\mathbb{P}^1$ & $L_2$ error        & EOA & $L_2$  error       & EOA & $L_2$ error       &EOA  \\ \hline 
                        20    & 5.98E-04  &  --  & 6.57E-04  &  --  & 6.89E-04  &  --    \\
                        40    & 1.27E-04  & 2.23 & 1.32E-04  & 2.31 & 1.34E-04  &  2.35  \\
                        80    & 3.02E-05  & 2.07 & 3.05E-05  & 2.10 & 3.07E-05  &  2.12  \\
                        160   & 7.45E-06  & 2.01 & 7.47E-06  & 2.03 & 7.48E-06  &  2.03  \\
                        320   & 1.86E-06  & 2.00 & 1.86E-06  & 2.00 & 1.86E-06  &  2.00  \\ \hline
        \end{tabular}
        \begin{tabular}{c|cc|cc|cc} 
                        &\multicolumn{2}{c|}{$d=-1$, CFL$=100$} &\multicolumn{2}{c|}{$d=-0.5$, CFL$=100$} &\multicolumn{2}{c}{$d=+0.8$, CFL$=100$}  \\[0.5mm]
                        \cline{2-7}
                        $N_e$, $\mathbb{P}^2$ & $L_2$ error        & EOA & $L_2$  error       & EOA & $L_2$ error       &EOA  \\ \hline 
                        20    & 2.68E-03 &  --  & 7.43E-04  &  --  & 5.89E-05  &  --   \\
                        40    & 3.37E-04 & 3.00 & 9.34E-05  & 3.00 & 7.43E-06  & 2.98  \\
                        80    & 4.22E-05 & 3.00 & 1.17E-05  & 3.00 & 9.31E-07  & 3.00  \\
                        160   & 5.28E-06 & 3.00 & 1.46E-06  & 3.00 & 1.16E-07  & 3.00  \\
                        320   & 6.60E-07 & 3.00 & 1.83E-07  & 3.00 & 1.45E-08  & 3.00  \\ \hline
        \end{tabular}
        \centering
        \begin{tabular}{c|cc|cc|cc} 
                        &\multicolumn{2}{c|}{$d=-1$, CFL$=100$} &\multicolumn{2}{c|}{$d=-0.5$, CFL$=100$} &\multicolumn{2}{c}{$d=+0.8$, CFL$=100$}  \\[0.5mm]
                        \cline{2-7}
                        $N_e$, $\mathbb{P}^3$ & $L_2$ error     & EOA & $L_2$  error       & EOA & $L_2$ error       &EOA  \\ \hline 
                        20    &  2.59E-05 &  --  & 6.62E-06 &  --  & 3.40E-07 &  --  \\
                        40    &  8.12E-07 & 4.99 & 2.08E-07 & 4.99 & 1.59E-08 & 4.41 \\
                        80    &  2.54E-08 & 5.00 & 6.53E-09 & 4.99 & 8.89E-10 & 4.16 \\
                        160   &  7.65E-10 & 5.04 & 2.01E-10 & 5.02 & 5.06E-11 & 4.13 \\
                        320   &  1.07E-11 & 6.16 & 5.00E-12 & 5.33 & 4.43E-12 & 3.51 \\ \hline
                        \hline
        \end{tabular}
\end{table}

\section{Conclusions}\label{sec:conclusions}

This paper presented a detailed stability analysis, in the context of linear hyperbolic equations, 
of high-order discontinuous Galerkin (DG) methods coupled with the shifted boundary polynomial correction, an embedded boundary treatment.
The study, conducted by examining the eigenvalue spectrum of the discretized operators for a simplified two-element system, 
has provided critical insights into the stability constraints of this high-order, geometrically unfitted approach.
The main findings show that, in the context of boundary conditions imposed at a distance of the order of $\Delta x$ from the mesh interface, 
the shifted boundary polynomial correction introduces notable stability restrictions for explicit DG schemes, 
and the stability region is not symmetric when considering external ($d<0$) or internal ($d>0$) embedded boundaries.
For $d<0$, stability is maintained but requires a progressively more restrictive CFL condition as $d\rightarrow-\Delta x$ and the polynomial order $p$ increases.
For $d>0$, the method becomes unconditionally unstable beyond a critical distance $d>d_{max}(p)$, where $d_{max}$ decreases with increasing polynomial order.
Implicit time integration effectively overcomes the stringent CFL constraints associated with explicit methods.
For $d<0$, the scheme is unconditionally stable, representing a major improvement.
For $d>0$, the instability threshold $d_{max}$ is significantly increased compared to the explicit case. 
Furthermore, unconditional stability for all $d$ can be recovered by using a CFL number larger than a certain threshold, which depends on the polynomial degree.
An extensive suite of numerical tests for the linear advection equation with a manufactured solution perfectly corroborated the theoretical stability analysis.
The predicted convergence rates were recovered for all stable configurations $(d,\text{CFL})$, and instabilities were observed precisely where the analysis predicted them.

Future work will extend this stability analysis to other computational frameworks, for instance for stabilized continuous Galerkin in the context of hyperbolic equations, 
which appear to be more challenging than elliptic problems.
Further implementations of the implicit unconditionally stable method will be carried out for linear and nonlinear multi-dimensional systems to assess the generalizability of this 
one-dimensional analysis to higher dimension.
The shape of the analytical eigenvalues \eqref{eq:eigen} show that the instabilities only arise from the boundary cell, therefore the implementation of an explicit-implicit method in the spirit of \cite{may2017explicit}, where only the boundary cells are treated implicitly, may greatly improves the computational performances of the scheme. 
The possibility of perfoming a stability analysis in the sense of the $L_2$-norm, or energy stability \cite{carpenter1999stable,mattsson2003boundary}, will also be investigated. The latter may give insights on which terms are responsible of the instabilities.


\section*{Acknowledgments}
The author wishes to express sincere gratitude to Sixtine Michel and Davide Torlo for their careful reading of the manuscript and for their valuable constructive feedback. The author also thanks Armando Coco, Wasilij Barsukow and Emanuele Macca for the stimulating discussions on the topic of boundary conditions.

\bibliographystyle{plain}
\bibliography{references}
\end{document}